\documentclass[10pt]{article} 
\usepackage[preprint]{neurips_2023}

\usepackage{verbatim}
\usepackage[english]{babel}
\usepackage[utf8x]{inputenc}
\usepackage[protrusion=true,expansion=true]{microtype}

\usepackage{placeins}  
\usepackage{xcolor}
\usepackage{paralist}
\usepackage{amsmath,amssymb, empheq} 

\usepackage[displaymath, mathlines,running]{lineno}

\usepackage{algorithm}
\usepackage{algpseudocode}
\usepackage{fontawesome5}
\usepackage{hyperref}
\usepackage{amsfonts,nccmath,bbm,mathrsfs,mathtools}
\usepackage{caption} 

\usepackage{caption}     

\usepackage{amsthm}

\makeatletter
\renewcommand{\@noticestring}{Preprint.}
\makeatother

\newtheorem{problem}{Problem}[] 

\usepackage{appendix}
\usepackage{etoolbox} 
\usepackage{tikz-cd}
\newcommand*\linenomathpatch[1]{%
  \cspreto{#1}{\linenomath}%
  \cspreto{#1*}{\linenomath}%
  \csappto{end#1}{\endlinenomath}%
  \csappto{end#1*}{\endlinenomath}%
}
\newcommand*\linenomathpatchAMS[1]{%
  \cspreto{#1}{\linenomathAMS}%
  \cspreto{#1*}{\linenomathAMS}%
  \csappto{end#1}{\endlinenomath}%
  \csappto{end#1*}{\endlinenomath}%
}

\expandafter\ifx\linenomath\linenomathWithnumbers
  \let\linenomathAMS\linenomathWithnumbers
  \patchcmd\linenomathAMS{\advance\postdisplaypenalty\linenopenalty}{}{}{}
\else
  \let\linenomathAMS\linenomathNonumbers
\fi

\linenomathpatch{equation}
\linenomathpatchAMS{gather}
\linenomathpatchAMS{multline}
\linenomathpatchAMS{align}
\linenomathpatchAMS{alignat}
\linenomathpatchAMS{flalign}

\newtheorem{thm}{Theorem}
\newtheorem{lem}[thm]{Lemma}
\newtheorem{cor}[thm]{Corollary}

\theoremstyle{definition} 
\newtheorem{definition}{Definition}

\let\today\relax
\makeatletter
\def\ps@pprintTitle{%
    \let\@oddhead\@empty
    \let\@evenhead\@empty
    \def\@oddfoot{\footnotesize\itshape
         {} \hfill\today}%
    \let\@evenfoot\@oddfoot
    }
\makeatother

\usepackage{float}  

\newcommand{\bsym}[1]{\boldsymbol{#1}}
\newcommand{\norm}[1]{\lVert#1\rVert}

\newcommand{\Diff}{\text{Diff}(M)}
\newcommand{\Diffs}{\text{Diff}^s(M)}

\newcommand{\bu}{\bsym{u}}
\newcommand{\dist}{\text{dist}}
\newcommand{\dt}{\Delta t}
\newcommand{\bv}{\bsym{v}}

\newcommand{\id}{\text{id}}

\definecolor{Bleu}{rgb}{0.3, 0.2, 1.0}

\title{Diffeomorphic Neural Operator Learning}

\author{
  Seth Taylor\thanks{Corresponding Author} \\
  Department of Mathematics and Statistics\\
  McGill University\\
  Montreal, QC, H3A 0B9, Canada \\
  \texttt{seth.taylor@mail.mcgill.ca} \\
  \And
  Alex Bihlo \\
  Department of Mathematics and Statistics\\
  Memorial University of Newfoundland\\
  St. John’s, NL, A1C 5S7, Canada \\
  \texttt{abihlo@mun.ca} \\
  \And
  Jean-Christophe Nave \\
  Department of Mathematics and Statistics\\
  McGill University\\
  Montreal, QC, H3A 0B9, Canada \\
  \texttt{jean-christophe.nave@mcgill.ca} \\
}

\begin{document}

\maketitle

\begin{abstract}
We present an operator learning approach for a class of evolution operators using a composition of a learned lift into the space of diffeomorphisms of the domain and the group action on the field space. In turn, this transforms the semigroup structure of the evolution operator into a corresponding group structure allowing time stepping be performed through composition on the space of diffeomorphisms rather than in the field space directly. This results in a number of structure-preserving properties related to preserving a relabelling symmetry of the dynamics as a hard constraint. We study the resolution properties of our approach, along with its connection to the techniques of diffeomorphic image registration. Numerical experiments on forecasting turbulent fluid dynamics are provided, demonstrating its conservative properties, non-diffusivity, and ability to capture anticipated statistical scaling relations at sub-grid scales. Our method provides an example of \emph{geometric operator learning} and indicates a clear performance benefit from leveraging a priori known infinite-dimensional geometric structure.
\end{abstract}


\section{Introduction}

The successes of deep learning  \cite{lecun2015deep} for image and natural language processing \cite{krizhevsky2012imagenet, brown2020language} have propelled its application to scientific computing.  By leveraging the expressivity of deep neural network approximations along with the availability of rich open-source data sets, these techniques present a powerful tool to forecast nonlinear dynamics. A validating example comes from the field of numerical weather prediction where training on ERA5 global re-analysis \cite{hersbach2020era5} and model outputs has produced forecasts which rival and excel the skill of models built upon established numerical methods \cite{bi2023accurate, kurth2023fourcastnet, lam2023learning}. \par

The development of operator learning methods are at the forefront of this progress, with the design of neural operator approximations receiving significant attention \cite{bhattacharya2021model, li2020fourier, lu2021learning, raonic2023convolutional}. In these methods, there has been a particular emphasis on the design of functional schemes \cite{bhattacharya2021model, li2020fourier, batlle2024kernel}, resulting in operators which can be evaluated on previously unseen data at varying resolution for downstream tasks such as multi-scale parameter optimization and uncertainty quantification. Developing techniques to constrain these functional approaches based on our extensive physical knowledge of the observed dynamics is fundamental to producing more reliable and accurate predictions \cite{kashinath2021physics}. This is especially imperative for the long-term prediction of global atmospheric dynamics, relevant to climate modelling, where it is essential to capture the delicate thermodynamic balances underlying the physics of the Earth system \cite{lauritzen2022reconciling}.  \par

\subsection{Structure-preserving operator learning} Hybrid modelling strategies which incorporate domain knowledge to physically constrain the learning problem can be broadly classified into indirect and direct approaches. Indirect techniques use pointwise constraints in the objective functional in the form of partial differential equations for instance \cite{beucler2019achieving, raissi2019physics, li2024physics}. These `soft' constraints however do not immediately generalize to unseen data and are sensitive to the training. Direct methods impose global constraints into the learning problem via the careful design of specialized architectures \cite{mohan2020embedding, kashinath2020enforcing} or through geometric approaches which embed symmetries and invariances \cite{wang2020incorporating, bronstein2017geometric}. In the setting of operator learning, direct methods require imposing an infinite-dimensional geometric structure or finite-dimensional analogue into the network architecture and/or data representation. Imposing structure is fundamental to producing computationally tractable algorithms in this setting; operator learning of general $C^k$ maps between Banach spaces are known to exhibit curses of parametric complexity \cite{lanthaler2023curse, novak2009approximation} and data complexity \cite{kovachki2024data}. Unlike the finite-dimensional setting, smoothness alone is insufficient to develop more accurate algorithms and it is only through additional structure \cite{adcock2022sparse, cohen2015approximation, marcati2023exponential}, or known representation formula \cite{lanthaler2022error}, that these curses of dimensionality manifest at the operator level have been shown to be overcome. \par

In their recent work, Bouziani and Boull{\'e} \cite{bouziani2024structure}, developed a structure-preserving operator learning approach which constrains the output of the model to a finite-dimensional conforming approximation space of functions. This approach directly embeds hard constraints into the learning problem, coming in the form of boundary conditions, domain geometry, and even discrete differential structures \cite{arnold2018finite} could be imposed directly through the use of known finite element discretizations. Accurate representations on a fixed mesh in a linear space of finite-elements are however challenged by the generation of fine scales in turbulent fluid dynamics and a variety of mesh-less Lagrangian techniques have been proposed to accommodate. 

\textbf{Lagrangian Method Emulation.} The quantity of interest at the core of most Lagrangian methods is the flow map associated to an Eulerian velocity field. Let $(M,g)$ be a $d
$-dimensional compact Riemannian manifold modelling the fluid domain and define the forward trajectory (Lagrangian) map $\varphi_{[0,t]}: M \to M$ through the ordinary differential equation
\begin{equation}
\label{flow_equation}
\partial_t \varphi_{[0,t]} = \bv(t) \circ \varphi_{[0,t]} \,, \quad \varphi_{[0,0]} = \id_M \,.
\end{equation}
where $\bv(t): M \to TM$ is a vector field on the domain. The forward Lagrangian map possesses an inverse $\varphi_{[t,0]}:M \to M$ such that $\varphi_{[t,0]}\circ\varphi_{[0,t]} = \id_M$ which can be expressed through an ordinary differential equation of the form
\begin{equation}
\label{flow_equation_2}
\partial_t \varphi_{[t,0]} = \hat{\bv}(t) \circ \varphi_{[t,0]} \,, \quad 
\hat{\bv}(t) =  -D\varphi_{[0,t]}^{-1} \cdot \bv(t) \circ \varphi_{[0,t]} \,,
\end{equation}
or through the transport equation
\begin{equation}
\label{map_transport}
\partial_t\varphi_{[t,0]} + D\varphi_{[t,0]}(\bv(t)) = 0  \, ,
\end{equation}
in the Eulerian frame where $D\varphi_{[t,0]}: TM \to TM$ is its differential. The inverse map is commonly of interest as it defines a solution operator to conservation laws of the form
\begin{equation}
\label{conservation_law}
\partial_t \rho + \nabla \cdot (\rho \bv) = 0\,, \quad \rho_0 \in L^2(M)\,,
\end{equation}
via $\rho_0 \mapsto  \rho_0 \circ \varphi_{[t,0]}\cdot J_{\mu}(\varphi_{[t,0]})$ with $J_{\mu}(\varphi_{[t,0]})$ being the Jacobian determinant of the inverse map defined by the metric volume form. In the cases where the underlying vector field is unknown, yet prior knowledge in the form \eqref{conservation_law} is given, then the backward and forward flow maps can be learned to directly enforce the integral conservation law associated with the evolution \eqref{conservation_law}. Existing Lagrangian numerical method emulation strategies have been devised on these discretizations such as the LFlow method \cite{torres2023lagrangian} using the inverse map approach along with diffusion models for Lagrangian turbulence \cite{li2024synthetic}. We refer to \cite{ wang2024recent} for a review on the state of the art for machine-learning in fluid mechanics. \par 

 \textbf{Nonlinear Approximation Space.} The method presented in this work offers a Lagrangian method emulation technique using composite approximation spaces for diffeomorphisms based on the recently developed Characteristic Mapping (CM) method \cite{yin2021characteristic, yin2023characteristic, taylor2023projection}. The method utilizes a spatiotemporal discretization of the flow map in a nonlinear composite approximation space; consider partitioning an interval of time $[0,t]$ into subdivisions $[\tau_{i}, \tau_{i+1}] \subset [0, t]$, we can decompose the paths formed by the inverse and forward maps into a composition of sub-interval flows of the form 
\begin{equation}
\begin{aligned}
\label{eq:global_map}
\varphi_{[0,\tau_n]} &= \varphi_{[\tau_{n-1}, \tau_n]} \circ \dots \circ \varphi_{[\tau_1, \tau_2]} \circ \varphi_{[0, \tau_1]}
\\
\varphi_{[\tau_n,0]} &= \varphi_{[\tau_1, 0]} \circ \varphi_{[\tau_2, \tau_1]}  \dots \circ \varphi_{[\tau_n, \tau_{n-1}]}\,,
\end{aligned}
\end{equation}
where the submaps solve the ODEs \eqref{flow_equation} and \eqref{flow_equation_2} over the subintervals respectively. The decomposition of the solution to \eqref{map_transport} yields a representation of a complex deformation of the domain as a composition of simpler and hence more accurately computed deformation maps. The flow maps define bijective mappings of the domain in the space of $C^1$ diffeomorphisms
\begin{equation}
\label{Diff}
C^1\Diff = \left\{\varphi \in C^1(M,M) \,|\, \varphi \text{ is a bijection}, \, \varphi^{-1} \in C^{1}(M,M)  \right\},
\end{equation}
which serves as the model approximation space. The space \eqref{Diff} can be treated as an infinite-dimensional Lie group \cite{schmid2004infinite, schmeding2022introduction} with the group action of composition and the flows of a velocity field form a one-parameter subgroup. Our approach uses a conforming approximation space of diffeomorphisms $D_h \subset C^1\Diff$ which serves as a building block to approximate complex deformations using the nonlinear approximation spaces
\begin{equation}
\label{composite_discretization}
D_h^{\circ, k} = \underbrace{D_h \circ D_h \circ \dots \circ D_h}_{k-\text{times}}  \subset C^1\Diff \,.
\end{equation}
Each layer forming the composite discretization \eqref{composite_discretization} results in an increasingly multi-scale approximation of the flow map governing the evolution. In turn, this accommodates for the generation of fine scales ubiquitous in the forecasting of nonlinear dynamics. The use of a conforming discretization of diffeomorphisms imposes hard constraints on the conservation principles associated with transported fields, stemming from a relabeling symmetry of the dynamics. Our proposed operator learning technique offers a structure-preserving approach capable of directly enforcing this infinite-dimensional geometric constraint.\par 

\subsection{Contributions}

The main contributions of this work are associated with the design of a novel structure-preserving operator learning strategy and its application for learning nonlinear fluid dynamics. We propose a geometric approach to evolution operator learning based on a lifting into the group of diffeomorphisms acting on the field space. This technique offers a number of unique resolution properties and performance implications for time-stepping which we elucidate. The time-stepping methods are immediately amenable to the use of existing architectures with a simple modification to output into an approximation space of maps of the domain. We validate our approach numerically and observe exceptional performance of the method in resolving turbulent fluid motions at subgrid scales without introducing diffusivity. The use of the nonlinear composite function space \eqref{composite_discretization} results in super-resolution properties which are different, yet complementary, to those observed in existing neural operator architectures \cite{kovachki2023neural}. We connect these properties to infinite-dimensional geometric structures preserved by our technique. 

\noindent
\textbf{Organization.} In Section \ref{sec:background_motivation} we begin by contextualizing our contribution and the methodological problems we seek to address. In Section \ref{sec:method} we then describe the method in detail, elaborate on its resolution properties, and outline a connection to the techniques of diffeomorphic image registration. We provide a numerical verification of these properties in section \ref{sec:numerical_verification} along with comparative numerical experiments of our approach against the Fourier Neural Operator architecture are given for learning the solution operator to the transport equation and the incompressible Euler equations. We conclude \ref{sec:conclusions} with a discussion on the limitations of the method along with future directions to the design of other geometric operator learning methods inspired by these techniques.  

\section{Background and Motivation}
\label{sec:background_motivation}

\subsection{Operator Approximations}

Let $U,V$ be two separable Banach spaces of functions $u : M \to \mathbb{R}^{m}$ on an $d$-dimensional compact Riemannian manifold $(M,g)$ such that $U \subseteq V$ continuously and denote the space of Lipschitz continuous operators as $C(U,V)$. The operator learning problem requires the definition of a parametric operator-valued map of the form
\begin{equation}
\label{functional_discretization}
\mathscr{D}: \Theta \subset \mathbb{R}^p \to C(U,V) \,,
\end{equation}
and a means of selecting an optimal parameter to approximate a given $\Phi \in C(U,V)$. In the data-driven setting, without a priori knowledge on the evolution operator, the predominant methodology is the use of stochastic gradient descent for regularized empirical risk optimization problem leveraging a deep neural network approximation space \cite{kovachki2024operator}. In the encoder-decoder approach, the map \eqref{functional_discretization} is constructed using an embedding into a space of continuous functions defined by $\mathcal{N}: \Theta \to \mathcal{A} \subset C(\mathbb{R}^{d_U}, \mathbb{R}^{d_V})$ along with pairs of continuous maps $\mathcal{E}: U \to \mathbb{R}^{d_U}$ and $\mathcal{E}^*:\mathbb{R}^{d_U} \to U$,  $\mathcal{D}: \mathbb{R}^{d_V} \to V$ and $\mathcal{D}^*: V \to \mathbb{R}^{d_V}$ such that the projections
\begin{equation}
\label{projection}
\pi_U = \mathcal{E}^* \circ \mathcal{E}: U \to U_h \subset U \,; \quad \pi_V = \mathcal{D}^* \circ \mathcal{D}: V \to V_h \subset V \,, 
\end{equation}
are approximations of the identity on $U$ and $V$ respectively. The encoder-decoder pair defines a map \eqref{functional_discretization} of the form
\begin{equation}
\label{encoder_decoder_approx}
\mathscr{D}_{E}: \Theta  \to C(U,V) \,,  \quad \theta \mapsto \mathcal{D} \circ \mathcal{N}(\theta) \circ \mathcal{E} =  \mathcal{D} \circ \mathcal{N}(\theta) \circ \mathcal{E}\,.
\end{equation}
Neural operators (NOs) offer an alternative approach based on pointwise analogues of a finite-dimensional neural network in the form
\begin{equation}
\label{neural_operator}
\mathscr{D}_{N}: \Theta \to C(U,V) \,, \quad \theta \mapsto \mathcal{Q} \circ \mathcal{L}_{L} \circ \mathcal{L}_{L-1} \circ \dots \circ \mathcal{L}_1 \circ \mathcal{R} 
\end{equation}
where $\mathcal{R}: U \to U'$ is a lifting operator into a Banach space of functions $v: M \to \mathbb{R}^{d'}$ with $d' > d$ and $Q: U_L' \to V$ is a projection operator. The layers $\mathcal{L}_k : U_k' \to U_{k+1}'$ of the neural operator are defined by
\begin{equation}
\label{neural_operator_layers}
\mathcal{L}_k[v](x) = \sigma(W_k\cdot v(x) + b_k(x)  + K_{\theta_k}[v](x))
\end{equation}
for all $x \in M$ where $W_k \in \mathbb{R}^{d'_{k} \times d_{k-1}'}$ and $b_k(x) \in \mathbb{R}^{d_k'}$ are the weights and biases, $\sigma\in C^{\infty}(\mathbb{R}^{d_k'}, \mathbb{R}^{d_k'})$ is the application of a nonlinear activation function component-wise, and the $K_{\theta_k}$ are parametric integral transforms, each performed at the $k$-th layer. The FNO \cite{li2020fourier} is a particular example of the neural operator \eqref{neural_operator} with an integral transform of the form
\begin{equation*}
\mathcal{K}_{\theta_k}[v](x) = \int_M \kappa_{\theta_k}(x - y)v(y) dy \,,
\end{equation*}
where $\kappa_{\theta}$ is a parametric integral kernel. Using the convolution identity of the Fourier transform $\mathcal{F}: L^2(\mathbb{T}^d, \mathbb{R}^m) \to \ell^2(\mathbb{Z}^d, \mathbb{C}^m)$ to apply the kernel as a multiplication in Fourier space this can be expressed as 
\begin{equation}
\label{fourier_operator}
\mathcal{K}_{\theta}[v](x) = \mathcal{F}^{-1}\left[\mathcal{F}[\kappa_{\theta}](k)\cdot \mathcal{F}[v](k)\right](x)  \,.
\end{equation}
In its implementation, the FNO requires a discrete Fourier transform, leading to a mixed approximation space which is neither purely local or purely spectral.

\subsection{Evolution Operator Learning}

The evolution operator learning problem considers the approximation of the solution operator $\Phi_t: U \to V$ to an initial value problem
\begin{equation}
\label{dynamical_system}
\dot{u}(t) = F(u(t)) \,, \quad u(0) = u_0\,,
\end{equation}
for an unknown $F \in C(U,V)$ which is locally Lipschitz continuous. The approximation is constructed from a supervised learning problem based on a finite sample of trajectories over an interval $[0,T]$ of time using a semi-discretized dataset
\begin{equation}
\label{dataset}
U_{N,\Delta t} = \{(u_i, \Phi_{k\Delta t}(u_i)) \in U \times V \,:\, 1  \leq i \leq N_d \,, \, 1 \leq k \leq N_t \} \,.
\end{equation}
\begin{problem}
\label{problem1}
 Construct a map \eqref{functional_discretization} which can approximate $\Phi_t : U \to V$ for some $t \in [0,T]$ given \eqref{dataset}.
\end{problem}
The problem \eqref{problem1} is substantiated by defining a notion of approximation, typically through a measure of the error over a distribution defined on the input space. In the construction of the map \eqref{functional_discretization}, the semi-group property of the evolution operator
\begin{equation}
\label{semi_group_property}
\begin{aligned}
\Phi_0(u) = u \,; \quad \Phi_{t+s}(u) = \Phi_t \circ \Phi_s (u) \,, \quad \forall u \in U \,,
\end{aligned}
\end{equation}
can be used to define an approximation of $\Phi_{\dt}:U \to V$ such that
\begin{equation}
\label{auto_regressive_time_stepping}
\Phi_t = \Phi^{\circ,n}_{\dt} \coloneqq \underbrace{\Phi_{\Delta t} \circ \Phi_{\Delta t} \circ \dots \circ \Phi_{\Delta t}}_{n-\text{times}} \,,
\end{equation}
with $t = n \dt$. The composite approximation \eqref{auto_regressive_time_stepping} distinguishes the evolution operator learning problem and the resulting spatiotemporal discretization defines different approximation spaces depending on the choice of functional discretization \eqref{functional_discretization}. In particular, we have that for all compact $K \subset U$ 
\begin{subequations}
\begin{align}
 \mathscr{D}_E(\theta) &= \tilde{\Phi}_{\dt} \implies  \tilde{\Phi}_{\dt}
^{\circ, n}(K) \subseteq V_h  \,, \label{encoder_decoder_composite}
\\
 \mathscr{D}_N(\theta) &= \tilde{\Phi}_{\dt} \implies  \tilde{\Phi}_{\dt}
^{\circ, n}(K) = \tilde{K}_n(\theta) \subset V \,, \label{neural_operator_composite}
\end{align}
\end{subequations}
where $\tilde{K}_n$ is a nonlinear subset which depends on the number of compositions taken and the approximation parameters. In the encoder-decoder approach, the control over the co-domain  \eqref{encoder_decoder_composite} offers better stability properties when aliasing is controlled \cite{mccabe2023towards, bartolucci2024representation, raonic2023convolutional} and it has even been demonstrated how learning on explicit discretizations of the domain and co-domain of the operator can result in better accuracy \cite{fanaskov2024spectral}. Longer time predictions however require larger band-limits to maintain stability imposing computational restrictions on the memory and efficiency. Moreover, de-aliasing the solution degrades the development of fine scale features and the nonlinear interactions from small to large scales. In contrast, the neural operator approximation results in a nonlinear approximation space \eqref{neural_operator_composite} which can be more expressive by adapting with the number of compositions taken. However, this increased expressivity often comes at the cost of reduced stability \cite{mccabe2023towards}, particularly when extrapolating over many time steps where the deviation from the training dataset accumulates. These properties of evolution operator learning highlight a need for schemes which are able to leverage the semigroup structure of the evolution operator while still balancing expressivity, stability, and computational efficiency across long time horizons.

\section{Diffeomorphic Neural Operator Learning}
\label{sec:method}

The main idea behind the proposed operator learning strategy is to \emph{learn how the field evolves under an action of the diffeomorphism group rather than how its values change}. The basic ansatz is to consider a continuous action of the diffeomorphism group on the field space   
\begin{equation}
\label{pullback_operator}
\Psi : \Diff \times U \to V \,, \quad (\varphi, u) \mapsto \varphi \cdot u \,,
\end{equation}
and a one-step parametric lifting operator
\begin{equation}
\label{lifting_operator}
\mathcal{K}_{\dt}:U \times \Theta \to C^1\Diff \,.
\end{equation}
This produces a functional discretization of the form 
\begin{equation}
\label{diffeomorphic_functional_discretization}
\mathscr{D}_{D}: \Theta \to C(U,V) \,, \quad \theta \mapsto \mathscr{C} \circ \mathcal{K}_{\dt}(\cdot, \theta)
\end{equation}
where we have defined the Lie group representation
\begin{equation}
\label{Lie_group_representation}
\mathscr{C}: C^1\Diff \to C(U,V) \,, \quad  \varphi \mapsto \Psi(\varphi, \cdot) \,.
\end{equation}
The construction \eqref{diffeomorphic_functional_discretization} for problem \eqref{problem1} exhibits a number of unique properties which can be related to two important features of the formulation. 1) The homomorphism property of the Lie group representation \eqref{Lie_group_representation} lifts the group structure of composition in $C^1\Diff$ into the space of operators allowing for the time-stepping to be performed completely in $C^1\Diff$. 2) Since the diffeomorphic group actions can be evaluated pointwise, the resulting operator approximations can act without reference to an explicit basis representation of the input function. This results in a consistency across resolution which is stronger than the discretization invariance exhibited by neural operators. Computationally, these properties are enabled by constraining the operator approximation to map into the nonlinear spaces defined by the orbits of the group action
\begin{equation}
\label{orbit_space}
\mathcal{O}_{u} = \left\{\Psi(\varphi,u) \,:\, \varphi \in C^1\Diff  \right\} \subset V \,,
\end{equation}
by construction. Invariants along the orbits can be associated to continuous relabelling symmetries of the dynamics and the functional discretization \eqref{functional_discretization} can be viewed as a structure-preserving approach similar to \cite{bouziani2024structure}. There is however there an essential distinction in that the inclusion $\iota_h: V_h \to V$ used in a conforming approximation space of finite-elements is replaced by the group action. Considering encoder-decoder approximation of the form $\mathcal{K}(\cdot, \theta) = \mathcal{D} \circ \mathcal{N}_{\theta} \circ \mathcal{E}$ where $\mathcal{N}_{\theta}: \mathbb{R}^{d_U} \to \mathbb{R}^{d_D}$ is a parametrized continuous function together with continuous maps of the form $\mathcal{E}: U \to \mathbb{R}^{d_U}$ and $\mathcal{D}:  \mathbb{R}^{d_D} \to D_h \subset C^1\Diff$, we can summarize this difference diagrammatically as:

\begin{equation}
\tikzset{ampersand replacement=\&}    
    \begin{tikzcd}
    \mathbb{R}^{d_U} \times U \arrow[r, "\tilde{\mathcal{K}}_{\theta}"] \& \mathbb{R}^{d_D} \times U \arrow[r, "{\mathcal{D}}"] \& D_h \times U \arrow[d, "\Psi"] \\%
    U \arrow[u, "\mathcal{E}"] \arrow{rr}{\Phi_t} \arrow[swap]{d}{\mathcal{E}} \& \& V  \\%
    \mathbb{R}^{d_U} \arrow{r}{\tilde{\Phi}_{\theta}} \& \mathbb{R}^{d_V} \arrow[r, "{\mathcal{D}}"] \& V_h \arrow[u, "\iota_h"']
    \end{tikzcd}
\end{equation}
This section elaborates on these resolution and time-stepping properties resulting from the approximation \eqref{diffeomorphic_functional_discretization} and establishes a connection for the construction of the lifting operator \eqref{lifting_operator} to the techniques of diffeomorphic image registration which provides a useful set of geometric and analytic tools to study the formulation.  

\subsection{Time Stepping via Group Action}
\label{sec:group_action_stepping}

Lifting a path in $\Diff$ into $C(U,V)$ through \eqref{Lie_group_representation} allows for the time-stepping to be performed by composing iterated lifts in place of the auto-regressive application \eqref{auto_regressive_time_stepping} directly on the field space. Using the shorthand $\mathcal{K}_{\dt}^{\theta}(u) = \mathcal{K}_{\dt}(u,\theta)$, denote 
\begin{equation}
\label{evolution_operator_decomposition}
\Phi_{\dt}(u) \approx \Psi(\mathcal{K}^{\theta}_{\dt}(u),u)
\end{equation}
for some fixed $\theta \in \Theta$. The semi-group property of the evolution operator \eqref{semi_group_property} and compatibility of the group action give distinct time-stepping schemes depending on how $\Diff$ acts on $U$:
\begin{equation*}
\begin{aligned}
\label{factorized_time_stepping}
\Phi^{\circ, n}_{\dt}(u)   \approx  \Psi(\varphi_n, \Psi(\varphi_{n-1},\cdots,\Psi(\varphi_1,u))) =  \begin{cases}
\Psi(\varphi_n \circ \varphi_{n-1} \circ \ldots \circ \varphi_1, u_0) \quad \text{if left} \,,
\\
 \Psi(\varphi_1 \circ  \ldots \circ \varphi_{n-1} \circ \varphi_n, u_0)  \quad \text{if right} \,.
\end{cases} 
\end{aligned}
\end{equation*}
where $\varphi_j = \mathcal{K}^{\theta}_{\dt}(\Psi(\varphi_{j-1}, u))$. Due to the direction of the composition in each case, an evaluation on $m$ points at every time step requires $\mathcal{O}(m)$ operations for the left action and $\mathcal{O}(n \cdot m)$ for the right action. Note that the composition operator $(\varphi,u) \mapsto u \circ \varphi$  is a right action whereas composition with the inverse $(\varphi, u) \mapsto u \circ \varphi^{-1}$ is a left action. The introduction of the inverse operation can however be costly and unstable for multiple compositions and we consider two schemes to perform the iteration \eqref{factorized_time_stepping} in the case of a right action. \par
The first scheme accumulates a history of the trajectory through the composition of maps, resulting in the time stepping
\begin{equation}
\label{composition_time_stepping}
\begin{aligned}
\tilde{\varphi}_{[t_{n+1},t_{n}]} &= \mathcal{K}^{\theta}_{\Delta t}(\tilde{u}(t_n))
\\
\tilde{u}(t_{n+1}) &= \Psi(\tilde{\varphi}_{[t_1, 0]} \circ \dots \circ \tilde{\varphi}_{[t_{n+1}, t_{n}]}  \,, u_0)
\end{aligned}
\end{equation}
initialized with $\tilde{u}(t_0) = u_0$. The scheme \eqref{composition_time_stepping} is a backward Lagrangian evolution technique which can be implemented in two variants. Either the trajectory is stored over the entire history of the time steps or the full composition chain is recomputed at each time step. The former results in an $\mathcal{O}(N_t \cdot |D_h|)$ memory allocation where $|D_h|$ is the number of degrees of freedom for each map, and the latter results in an $\mathcal{O}(N_t^2 \cdot N_x)$ computational cost to evolve $N_t$ time steps at $N_x$ grid points. \par

In cases where the full composition-based scheme \eqref{composition_time_stepping} becomes too memory intensive or computationally expensive over long time horizons, a semi-Lagrangian variant can be introduced to reduce storage and computation. This approach introduces an intermediate projection $\Pi_h: C^1 \Diff \to D_h \subset C^1\Diff$ onto an approximation space of diffeomorphisms at a subsequence of the time steps $\{\tau_j\} \subset \{t_i\}$ which can be chosen adaptively or a priori. Suppose that each macro time interval $[\tau_j, \tau_{j+1}]$ is subdivided by time steps $\tau_{j} =  t^j_1 < \dots < t_{m(j)}^j = \tau_{j+1}$. The semi-Lagrangian time stepping scheme can then be written as 
\begin{equation}
\label{SL_composite_time_stepping}
\begin{aligned}
\tilde{\varphi}_{[\tau_{j}, \tau_{j-1}]} & = \Pi_h[\mathcal{K}^{\theta}_{\dt}(u(t^j_{1})) \circ \mathcal{K}^{\theta}_{\dt}(u(t^j_{2})) \circ \dots \circ \mathcal{K}^{\theta}_{\dt}(u(t^j_{m(j)}))]
\\
\tilde{\varphi}_{[t_{n+1},\tau_{j}]} &= \tilde{\varphi}_{[t_n, \tau_j]} \circ \mathcal{K}^{\theta}_{\Delta t}(u(t_n))
\\
u(t_{n+1}) &= \Psi(\tilde{\varphi}_{[\tau_1, 0]} \circ \tilde{\varphi}_{[\tau_2, \tau_1]} \circ \dots \circ \tilde{\varphi}_{[\tau_{j},\tau_{j-1}]} \circ \tilde{\varphi}_{[t_{n+1},\tau_j]}, u_0)  \,.
\end{aligned}
\end{equation}
where the intermediate maps are initialized as $\varphi_{[\tau_j, \tau_j]} = \id$ after the projection step is taken. 
  
\subsection{Resolution Properties}
\label{sec:resolution_properties}

In this section we elaborate on some of the unique resolution properties resulting from the construction of the functional discretization \eqref{diffeomorphic_functional_discretization} and the nonlinear approximation spaces
\begin{equation}
\label{discrete_orbit_space_constraint}
\mathcal{O}_{u_0}^{\circ, k} \coloneqq \left\{\Psi(\varphi, u_0) \,: \, \varphi \in D^{\circ,k}_h \,\right\} \subset \mathcal{O}_{u_0} \,,
\end{equation}
defined by the schemes of section \ref{sec:group_action_stepping}. 

\textbf{Composite Discretization.} The discrete orbit spaces \eqref{discrete_orbit_space_constraint} possesses an exponentially growing bandwidth in the number of compositions used in the nonlinear space \eqref{composite_discretization}. We can coarsely characterize this property by considering the growth of the support of in frequency space as a function of the number of compositions. Let $M = \mathbb{T}^d$ and denote the space of band-limited functions as 
\begin{equation}
\label{band_limited_functions}
\mathcal{B}_{L} = \left\{f \in L^2(\mathbb{T}^d) \,: \text{supp}(\mathcal{F}(f)) \subseteq [-L,L] \right\} \,.
\end{equation}
The diffeomorphism group $C^1\text{Diff}(\mathbb T^d)$ can be identified with the space 
\begin{equation}
\label{approx_diff_space}
\mathfrak{X}^1_d(\mathbb{T}^d) = \left\{\varphi : \mathbb{T}^d \to \mathbb{R}^d \,:\, \varphi - \id \in  C^1(\mathbb{T}^d, \mathbb{R}^d) ,\,\,\,  \text{det}(D\varphi) > 0\right\}
\end{equation}
through the map 
\begin{equation}
\label{quotient_space_chart}
J : \mathfrak{X}^1_d(\mathbb T^d) \to C^1\Diff \,, \quad \bv \mapsto \pi(\id + \bv)
\end{equation}
where $\pi: \mathbb{R}^d \to \mathbb{T}^d \simeq \mathbb{R}^d/2\pi \mathbb Z^d$ is the pointwise$\mod 2 \pi$ operation in each coordinate direction. Define the $\epsilon$-effective bandwidth as
\begin{equation}
\label{epsilon_bandwidth}
\text{bw}_{\epsilon}(u) = \min\bigg\{R \in \mathbb{N}\,:\, \sum_{|k| \geq R} |\widehat{u}|^2 \leq \epsilon^2 \bigg\} \,.
\end{equation}
\begin{lem}
\label{lem:composite_bandwidth}
Approximating \eqref{approx_diff_space} using \eqref{band_limited_functions} in each coordinate direction gives 
$D_h \subset \pi(\id + \mathcal{B}_L^d)$ with
\begin{equation}
\label{composite_bandwidth}
\text{bw}_{\epsilon}(\varphi - \id) \le\; L \left[1 + \log(C_1/\epsilon) \right] \frac{(1 + C_2L)^k -1}{LC_2} \,, \quad \forall \varphi \in D_h^{\circ,k} \,,
\end{equation}
where $C_1 = \max_{1 \leq j \leq k}\norm{\varphi_j - \id}_2$ and $C_2 = \max_{2 \leq j \leq k}\norm{\varphi_j \circ \varphi_{j+1} \circ \cdots \circ \varphi_k - \id}_{\infty}$, and   
\begin{equation}
\label{orbit_space_estimate}
\text{bw}_{\epsilon}(u) \leq L_0 \left[1 + \log(C_1/\epsilon) \right] \frac{(1 + C_2L)^k -1}{LC_2} \,, \quad \forall u_0 \in \mathcal{B}_{L_0} \,, u \in \mathcal{O}_{u_0}^{\circ,k} \,.
\end{equation}
\end{lem}
The proof is provided in the Appendix \ref{app:proofs}. Lemma \ref{lem:composite_bandwidth} indicates how the composite discretization spaces \eqref{composite_discretization} possess an exponentially increasing bandwidth with only a linear increase in degrees of freedom since $|D_h^{\circ,k}| = k |D_h|$. This results in an exponentially increasing expressivity using only a linear increase in computational resources.

\textbf{Resolution Consistency.} The continuum approximation used in the neural operator formulation \cite{li2020fourier, kovachki2023neural, kovachki2024operator} provides training strategies which can decouple the resolution of the dataset from the parametrization of the operator. This has been shown to exhibit improved accuracy at inference using higher resolution data \cite{wei2023super, li2020fourier} and framed using the notion of discretization invariance based on a discrete refinement \cite{kovachki2023neural}. A \emph{discrete refinement} of the domain $M$ is a sequence of nested sets $M_1 \subset M_2 \subset \dots \subset M$ with $|M_L| = L$ for any $L \in \mathbb{N}$ such that for any $\varepsilon > 0$, there exists a number $L(\epsilon) \in \mathbb{N}$ with
\begin{equation}
M \subseteq \bigcup_{x \in M_{L(\epsilon)}} \left\{ y \in M \,:\, d_g(x,y) < \epsilon \right\}\,.
\end{equation} 
Define the associated sampling map which evaluates a function on $M_L$ as $\mathscr{S}_L: V \to \mathbb{R}^{L \cdot m}$. 
Given a discrete refinement $\{M_L\}_{L = 1}^{\infty}$, a parametric operator approximation $\tilde{\Phi}: \Theta \times U \to V $ is called \emph{discretization invariant} if there exists a sequence of maps $\tilde{\Phi}_{L}: \mathbb{R}^{L \cdot d}\times \mathbb{R}^{L \cdot m} \times \Theta \to V$ such that for any $\theta \in \Theta$ and any compact set $K \subset U$ the following holds 
\begin{equation}
\label{discretization_invariance}
\lim_{L \to \infty} \sup_{u \in K} \norm{\tilde{\Phi}_L(M_L, \mathscr{S}_L(u), \theta) - \tilde{\Phi}(\theta, u)}_V = 0 \,.
\end{equation}
Discretization invariance implies convergence under mesh refinement to an operator in the continuum. It however does not imply a consistency across resolutions \cite{gao2025discretization}; although the neural operator defines a map $\mathscr{D}_N: \Theta \to C(U,V)$, this representation is realized only through its evaluation on a fixed discretization of the inputs. The limit \eqref{discretization_invariance} can be seen as the pointwise definition of the neural operator and its representation at each resolution can correspond to maps which do not coincide with the restriction of a single operator. Our diffeomorphic neural operator approach possesses a sufficient condition for discretization invariance which we refer to as resolution consistency. Resolution consistency is however stronger notion than discretization invariance and implies a convergence under discrete refinement in a way that its the approximated operators' outputs can be viewed as an evaluation of the same function at finer resolutions. We use the following definition. \par

\begin{definition}
\label{def:resolution_consistent}
Given a discrete refinement $\{M_L\}_{L = 1}^{\infty}$, a discretization invariant operator $\tilde{\Phi}:\Theta \times U \to V$ is called \emph{resolution consistent}, if the sequence $\{\tilde{\Phi}_L\}$ satisfying \eqref{discretization_invariance} is such that
\begin{equation}
\mathscr{S}_{L'}\left(\tilde{\Phi}_L(M_L, \mathscr{S}_L(u), \theta)\right) = \mathscr{S}_{L'}\left(\tilde{\Phi}_{L'}(M_{L'}, \mathscr{S}_{L'}(u), \theta)\right)\,,
\end{equation}
for all $L$ and $L' \leq L$. 
\end{definition}

\begin{lem}
\label{lem:resolution_consistency_Diff}
Any continuous map $F: \Theta \to D_h \subset C^1\Diff$ defines a resolution consistent operator 
\begin{equation}
\label{basic_diff_operator}
\tilde{\Phi} : \Theta \times U \to V \,, \quad (\theta, u) \mapsto \Psi(F(\theta), u) \,.
\end{equation}
\end{lem} 
The proof of Lemma \ref{lem:resolution_consistency_Diff} is straightforward and given in the Appendix \ref{app:proofs}. This relates to the construction \eqref{diffeomorphic_functional_discretization} with the main observation being that no additional parametrization or discretization error is introduced when we evaluate \eqref{basic_diff_operator}: once the element in $\Diff$ has been computed, the action can be applied
\emph{exactly} at any resolution. This is made possible by the fact that the output of the lifting operator \eqref{lifting_operator} \emph{is} an operator through the representation \eqref{Lie_group_representation}. In particular, for any fixed $u \in U$ the map $\Psi(\mathcal{K}(u,\cdot),\cdot):U \times \Theta \to V$ is resolution consistent. Using a neural operator approximation of the lifting operator \eqref{lifting_operator} yields a discretization invariant approximation which possesses two distinct ways of upsampling: either evaluating the output diffeomorphism on a finer grid and composing at these points or leveraging continuum formulation of the neural operator to generate higher resolution diffeomorphisms. The latter outputs operators which can be evaluated arbitrarily throughout the domain whereas the former is restricted to uniform samplings.

\textbf{Relabelling Symmetry and Conservation}. The orbit space constraint \eqref{discrete_orbit_space_constraint} relates to a relabelling symmetry of the dynamics. Elements of $\mathcal{O}_u$ inherit the properties of $u$ invariant under coordinate transformation. Let $\varphi(t):M \to M$ be a Lagrangian flow map and denote $\varphi(S_0) = S(t) \subset M$. Consider a conservation law in the form 
\begin{equation}
\label{conservation_law}
\frac{d}{dt} \int_{S(t)} u(t) = 0 \,,
\end{equation}
where $u(t) \in \Omega^k(M)$. Changing variables with $\tilde{\varphi} \in C^1\Diff$ gives 
\begin{equation}
\int_{\varphi(S_0)} \varphi_*u_0 = \int_{S_0} u_0 = \int_{\tilde{\varphi}(S_0)}\tilde{\varphi}_*u_0 \,,
\end{equation}
where $\varphi_* \coloneqq (\varphi^{-1})^*$. The integral over the initial condition is preserved along any diffeomorphism of the domain, however pointwise there is an \emph{error in location} of the transported set. This leads to an advective notion of error in the solution \cite{taylor2023characteristic, yin2023characteristic}: for any approximate $\tilde{u} = \tilde{\varphi}_*u_0 \approx \varphi_*u_0 = u$ there exists an error map $\varepsilon = \varphi \circ \tilde{\varphi} \in C^1\Diff$ such that $\varepsilon^*\tilde{u}(t) = u(t)$. The error map describes a bijective correspondence between the solution and the approximation. Consequently, the numerical approximation is non-diffusive, preserves the minimum, maximum, and critical points and possesses an analytically exact conservation for any globally conserved quantity such as mass. 

In the fluid dynamics community, material invariance of Lagrangian methods is a well-known property. Discrete conservation follows from computing the evolution along material point discretizations of the evolution with arbitrarily labelled data. These methods possess discrete relabelling symmetries defined by permutations of the label indices or more generally for other discrete Lagrangian elements \cite{bowman2015fully}. The mapping-based approach we use enforces instead a continuous relabelling symmetry owing to the use of the conforming discretization $D_h \subset C^1\Diff$.   

\subsection{Connection to Diffeomorphic Registration Methods} 
\label{sec:LDM_connection}
In order to characterize the quality of the approximation \eqref{diffeomorphic_functional_discretization} we must understand how to construct a decomposition of the form \eqref{evolution_operator_decomposition} with some theoretically `closest` lifting operator. This can be accomplished using a two-point boundary value problem formulation on the diffeomorphism group:  given $u \in U$ find a $\varphi \in \Diff$ which minimizes an energy functional of the form 
\begin{equation}
\label{energy_Diff}
E(\varphi) = \norm{\Phi_{\dt}(u) - \Psi(\varphi,u)}_U^2 + \mathcal{R}(\varphi) \,,
\end{equation}
where $\mathcal{R}:\Diff \to [0,\infty)$ is a regularization term used to ensure well-posedness. Variational problems involving energies of the form \eqref{energy_Diff} arise in a variety of nonlinear statistical problems for medical imaging \cite{ceritoglu2013computational, gerber2009manifold}, and a number of algorithms have been proposed for their solution \cite{trouve1998diffeomorphisms, dupuis1998variational, miller2006geodesic, ashburner2011diffeomorphic, vialard2012diffeomorphic}. These methods provide techniques for nonlinear statistical analyses between shapes \cite{younes2010shapes} and possess a rich geometric structure \cite{bauer2014overview} which can be leveraged for constructive approximation and analysis. 

Letting the energy functional \eqref{energy_Diff} be variable in $u \in U$, we can define a lifting operator in the form 
\begin{equation}
\label{lifting_operator_definition}
\mathcal{K}_{\dt}(u) =  \arg\min_{\varphi \in \Diff}  \,\, \norm{\Phi_{\dt}(u) - \Psi(\varphi,u)}_U^2 + \mathcal{R}(\varphi) \,.
\end{equation}
In order to ensure the well-posedness of this formulation there are deliberate functional analytic and topological considerations that must be taken for the construction of the regularization and with respect to the continuity of the map \eqref{lifting_operator_definition}. 

\textbf{Flows of Diffeomorphisms.}  In the flows of diffeomorphisms model \cite{beg2005computing}, a Riemannian metric on $\Diff$ is constructed to regularize the energy functional \eqref{energy_Diff} to the one-parameter subgroups defined by the flows of time-dependent vector fields. Let $\mathfrak{X}(M)$ denote the space of vector fields and note that the tangent spaces are the right-translation of vector fields, that is
\begin{equation*}
\label{tangent_space}
T_{\varphi}\Diff = \left\{\bv : M \to TM \,: \, \pi_{TM} \circ \bv = \varphi\right\} = \mathfrak{X}(M) \circ \varphi \,,
\end{equation*}
where $\pi_{TM}:TM \to M$ is the tangent bundle projection.  An inner product on the space of vector fields can then be extended to the rest of the tangent bundle to define an invariant metric of the form
\begin{equation}
\label{right_invariant_metric}
G_{\varphi}(\cdot, \cdot):T_{\varphi}\Diff \times T_{\varphi}\Diff \to \mathbb{R}  \,, \quad (\bu \circ \varphi, \bv \circ \varphi) \mapsto \int_M g(\bu, L\bv) \mu
\end{equation}
where $L:\mathfrak{X}(M) \to \mathfrak{X}(M)$ is a positive, elliptic, and self-adjoint differential operator which enforces some desired level of smoothness. Considering $L = (1 - \alpha \Delta)^{s}$ or $L = I - \Delta^{s}$ where $\Delta$ is the vector Laplacian, define the space
\begin{equation}
\label{admissible_vector_fields}
X^s([0,T]) = \left\{\bv :[0,T] \to \mathfrak{X}^s(M) \,:\, \int_0^T \norm{\bv(t)}_s \,dt< \infty   \right\} \,,
\end{equation}
where $\norm{\cdot}_s$ is the norm induced by the inner product and $\mathfrak{X}^s(M)$ is the closure of $\mathfrak{X}(M)$ with respect to this norm. Denote $\Diffs$ as the closure of $C^1\Diff$ in the $H^s(M,M)$ Sobolev metric topology with $s > d/2 + 1$. The metric \eqref{right_invariant_metric} induces a distance on $\Diffs$ in the form
\begin{equation}
\label{Hs_distance}
\dist_{s}(\varphi, \eta) = \underset{\bv \in X^s([0,T])}{\text{inf}} \left\{ \int_0^T \norm{\bv(t)}_s \, dt \,:\,  \varphi = \gamma(T) \,, \, \gamma(0) = \eta \,, \, \dot{\gamma}(t) = \bv \circ \gamma(t)  \right\} \,,
\end{equation}
which can be used to define the regularization $R(\varphi) = \dist_s(\id, \varphi)$ penalizing the distance from the identity.  Bruveris and Vialard \cite{bruveris2017completeness} showed that $(\Diffs, \dist_s)$ is a complete metric space which provides a number of useful properties to study the formulation \eqref{lifting_operator_definition}. \par

\textbf{Theoretical Considerations.}
Existence of a minimizer in \eqref{lifting_operator_definition} can be established using the direct method of the calculus of variations (\cite{younes2010shapes} Theorem 10.2). Uniqueness is however not guaranteed in general due, in part, to a non-trivial stabilizer subgroup
\begin{equation}
\label{stabilizer_subgroup}
\text{Diff}^s_u(M) = \left\{ \varphi \in \Diffs \,:\, \Psi(\varphi,u) = u \right\} \,.
\end{equation}
Instead we can consider the lifting operator \eqref{lifting_operator_definition} to be set-valued $\mathcal{K}_{\dt}: U \rightrightarrows  \Diffs$ and seek a continuous selector of the form
\begin{equation}
\label{continuous_rep_lift}
 \tilde{\mathcal{K}}_{\dt}: 
\tilde{U} \subset U \to \Diffs \,, \quad \text{s.t.   }\,\,\,\,  \tilde{\mathcal{K}}_{\dt}(u) \in \mathcal{K}_{\dt}(u) \,, \quad \forall u \in \tilde U \,.
\end{equation}
A globally continuous selector cannot be constructed if the stabilizers \eqref{stabilizer_subgroup} are not all conjugate. Even the construction of a locally continuous selector is complicated by the topology of the orbits, requiring a locally transverse slice so that each equivalence class has a unique representative which depends continuously on input field. Instead, we can relax the continuity condition for the selector and consider a weaker space in which the approximation will hold. Let $\mu \in \mathcal{P}(U)$ be a probability measure over $U$ and denote the $L^2(\mu)$ norm as
\begin{equation}
\label{L2_op_norm}
\norm{\Phi}_{L^2(\mu, U)}^2 = \int_U \norm{\Phi(u)}_V^2 d\mu(u) \,.
\end{equation}
This setting is considered by Lanthaler et al \cite{lanthaler2022error} in their analysis of DeepONets for operator learning and for the FNO by Kovachki et al. (Theorem 18 \cite{kovachki2021universal}). Measuring the approximation error with \eqref{L2_op_norm} allows us to consider discontinuous approximations for the selector \eqref{continuous_rep_lift} which need only be measurable. In order to make our discussion more concrete, we will consider the action of composition
\begin{equation}
\label{composition_operation_Hs}
\Psi: \text{Diff}^{s}(M) \times H^{s+1}(M) \to H^s(M) \,, \quad (\varphi, f) \mapsto f \circ \varphi \,,
\end{equation}
with $s > d/2 + 1$ to ensure differentiability \cite{inci2013regularity}. 
It will be convenient to work within a local chart on $\Diffs$ defined by a $C^{1}$ diffeomorphism of the form
\begin{equation}
J:\mathcal U \subset \mathfrak{X}^s(M) \to \tilde{\mathcal{U}} \subset \Diffs, \qquad J(0)=\id,
\end{equation}
with inverse $J^{-1}$ defined on $\mathcal U$ and on bounded subsets both $J$ and $J^{-1}$ are locally Lipschitz. Riemannian normal coordinates can be considered using the Riemannian exponential map defined by the metric of the form \eqref{right_invariant_metric} which is a local diffeomorphism on a neighbourhood of the identity \cite{misiolek2010fredholm}. In the case $M=\mathbb T^{d}$, we can define a global coordinate chart on the diffeomorphism group using \eqref{quotient_space_chart}. Note that the step size $\dt$ controls the distance from the identity since
\begin{equation}
\label{set_control_Diff}
\frac{1}{2}\dist^2_s(\id, \varphi) \leq E(\varphi;u) \leq E(\id; u) = \norm{\Phi_{\dt}(u) - u}_V^2 \leq \rho (\Delta t) \,,
\end{equation}
for any $\varphi \in \mathcal{K}_{\dt}(u)$ and $u \in U$ by the minimizing property, which allows us to impose bounds on the size of the sets satisfying \eqref{lifting_operator_definition} in the metric space $\Diffs$. Using this construction we get the following theoretical results, with proofs given in the Appendix \ref{app:proofs}.
\begin{lem}
\label{lem:measurable_selector}
There exists a measurable selector $\tilde{\mathcal{K}}_{\dt}: U \to \Diffs$ such that $\tilde{\mathcal{K}}_{\dt}(u) \in \mathcal{K}_{\dt}(u)$ for $\mu$-almost every $u \in U$.
\end{lem}
\begin{thm}
\label{thm:approx_thm1}
Let $\Phi_{\dt}: H^{s+1}(M) \to H^{s}(M)$ be a continuous operator and let $K \subset U$ be a bounded set. Then for any $\varepsilon > 0$ there exists a neural operator approximation $\mathcal{N}_{\theta}: K \to \mathcal{U} \subset \mathfrak{X}^s(M)$ satisfying 
\begin{equation}
\label{operator_approx_error}
\norm{\Phi_{\dt} - \Psi(J \circ \mathcal{N}_{\theta}(\cdot), \cdot)}_{L^2(\mu,K)}< \norm{\Phi_{\dt}(u) - \Psi( \tilde{\mathcal{K}}_{\dt}(\cdot), \cdot)}_{L^2(\mu,K)}  + \varepsilon
\end{equation}
\end{thm}
The first term in \eqref{operator_approx_error} reflects the model error, being the deviation from the group orbits, and the second is the approximation error. In the case that the orbits of the evolution operator are contained in the orbits of the group action, we get the following corollary.
\begin{cor}
\label{orbit_time_stepping_corollary}
Let $\Phi_{\dt}: H^{s+1}(M) \to H^s(M)$ be a continuous operator such that $\Phi_t(u) \in \mathcal{O}_u$ for all $t \in [0,T]$ and $u \in K \subset U$ with $K$ bounded and suppose that $\Phi_t(K) \coloneqq \bar{K} \subset B_{R}(0)$ for all $t \in [0,T]$. Then for any $\varepsilon > 0$ there exists a neural operator $\mathcal{N}_{\theta}: \bar{K} \to \mathcal{U} \subset  \mathfrak{X}^s(M)$ such that 
\begin{equation}
\int_{\bar{K}}\norm{\Phi_{\dt}(u) - \Psi \circ J \circ \mathcal{N}_{\theta}^{\circ,n}(u)}^2_{\infty} d\mu(u) < \varepsilon \,.
\end{equation}
and the accumulated $n$-step error is bounded by
\begin{equation}
\int_{\bar{K}}\norm{\Phi^{\circ,n}_{\dt}(u) - \Psi \circ J \circ \mathcal{N}_{\theta}^{\circ,n}(u)}^2_{\infty} d\mu(u) <  C\varepsilon\sum_{j = 0}^{n-1} (1 + \tilde{L})^{n-j-1} \,.
\end{equation}
where $\Psi \circ J \circ \mathcal{N}_{\theta}^{\circ,n}$ is the $n$-times application of $u_{k+1} = \Psi(J \circ \mathcal{N}_{\theta}(u_k),u_k)$ with $u_0 = u$ and where $\tilde{L}$ depends on the Lipschitz constants of $J$ and $\mathcal{N}_{\theta}$. 
\end{cor}
We note that by considering the pullback action on the space of normalized volume forms -- which is notably transitive \cite{moser1965volume} -- can simplify some of these technical difficulties and extensions to other diffeomorphic group actions is warranted. Nevertheless, these theoretical considerations highlight a need for locally defined operator approximations in this geometric setup. This is not unique to the use of the functional discretization \eqref{diffeomorphic_functional_discretization} but rather more fundamental to the geometry used in the construction. Further analysis and theoretical development is needed to establish the use of the tools of differential calculus in the construction of nonlinear approximations for operator learning.

\section{Numerical Experiments}
\label{sec:numerical_verification}

We provide two numerical experiments designed to illustrate the unique properties of the operator approximation \eqref{diffeomorphic_functional_discretization} considering the case $M = \mathbb{T}^2$ for simplicity. The first experiment demonstrates the method's ability to generalize to out of distribution discontinuous fields for data generated by the transport equation. The second experiment assesses the method's ability to learn turbulent fluid dynamics for data generated by the two-dimensional incompressible Euler equations and showcases the capacity to resolve the anticipated energy cascade scalings at subgrid scales. Each experiment clearly demonstrates the benefits of the hard constraints imposed by the method with challenging tasks for existing approaches due to their smoothing properties. 

\subsection{Lifting Operator Approximation}

We consider the use of a neural operator $\mathcal{N}: U \times \Theta\to \mathfrak{X}^1_d(\mathbb{T}^2)$ to approximate the one-step lifting operation in the form
\begin{equation}
\label{neural_lifting_operator}
\tilde{\mathcal{K}}_{\dt}: \Theta \times U \to  C^1\text{Diff}(\mathbb T^2) \,, \quad (\theta, u) \mapsto \pi \circ( \mathcal{N}(u,\theta) + \id) \,.
\end{equation}
In practice, the neural operator is evaluated on finite input data and the data used to train the neural operator is presented in the form \eqref{dataset} which is sampled on a uniform periodic mesh. The neural operator is constructed to output into a finite-dimensional subset of diffeomorphisms based on the uniformly sampled inputs. As approximation space we consider a space of vectorized Hermite cubic splines for the displacement component. Spline interpolation spaces are natural candidate approximation spaces for diffeomorphism approximation as they possess a locality of the basis functions, ensures that their composition can be performed efficiently, along with global differentiability necessary for the conforming property in $C^1\Diff$.  Define $\Omega_L$ as the uniform quadrilateral partitioning of the domain and define the approximation space as   
\begin{equation*}
D_h = \pi\left(\left\{f \in C^1(\mathbb{T}^2,\mathbb{R}^2) \,:\, \left.f\right\vert_{C_{i,j}} \in \mathbb{P}^2_3 \quad \forall C_{i,j} \subset \Omega \right\} + \text{id}_{\mathbb{T}^2}\right)\,.
\end{equation*}
The Hermite interpolation operator interpolant is constructed through a collection of linear functionals
\begin{equation}
\label{diff_interpolant}
\Pi_h^*: C^1(\mathbb{T}^2, \mathbb{R}^2) \to \mathbb{R}^{8\cdot N_v} \,, \quad F \mapsto \left\{(F(x_i), \partial_xF(x_i), \partial_y F(x_i), \partial_x \partial_yF(x_i)\right\}_{i = 1}^{N_v} \,,
\end{equation}
The dimensionality $N_v = N_v(L)$, defining the output mesh parameter $h = h(L)$ depends on the sampling parameter. The neural operator outputs onto $\mathbb{R}^{d \cdot N_v}$ given the input samples on a mesh defined by $L$. It is worth emphasizing that $\theta$ is not a function of $L$ however the training and thus the parametrization of the neural operator is a function of the resolution of the training data. The neural operator can be evaluated at different discretizations of the inputs after training has been performed.  The projection onto the approximation space is defined by  
\begin{equation}
 \Pi_h :  C^1(\mathbb{T}^2, \mathbb{R}^2) \to D_{h} \,, \quad  
F  \mapsto  \pi\left(\id + \sum_{i,j} \sum_{\alpha, \beta} F^{(\alpha,\beta)}_{i,j} H_{i,j}^{\alpha, \beta} \right)\,,
\end{equation}
where the basis functions are given by the tensor product
\begin{equation*}
H_{i,j}^{\alpha, \beta}(\bsym{x}) = Q_{\alpha}\left(\frac{x - x_i}{\Delta x}\right)Q_{\beta}\left(\frac{y - y_j}{\Delta y}\right)\Delta x^{\alpha + \beta}  \,, 
\end{equation*}
and the $Q$ are Hermite cubic polynomials (see \cite{yin2021characteristic} for more details). The interpolant is $\mathcal{O}(\Delta x^4)$ accurate in the $L^{\infty}$ norm with at least $C^4$ regularity and is globally differentiable.

\subsection{Experimental Setup and Implementation Details} 
The dataset is defined in the form \eqref{dataset} and we consider learning a one-step approximation of the evolution operator. This is performed by computing an optimal lifting neural operator \eqref{neural_lifting_operator} through a numerical optimization of the cost functional
\begin{equation*}
\begin{aligned}
\mathcal{C}(\theta) &= \frac{1}{N N_tN_x}\sum_{i = 1}^N\sum_{k = 1}^{N_t} \norm{u_{i}(t_k) \circ \mathcal{K}_{\dt}(u_i(t_k),\theta) - \Phi_{\dt}(u_i(t_k))}_{\ell^2} +  \lambda \sum_{\alpha}\norm{\partial^{\alpha}\mathcal{N}(u_i(t_k), \theta)}_{\ell^2} \,.
\end{aligned}
\end{equation*}
where $\alpha \in \{(0,1), (1,0), (1,1)\}$, and $\lambda >0$ is a tuning constant for the regularization. This penalizes only the derivative components of the output diffeomorphisms at the grid points defining the interpolants, ensuring the transformations remain sufficiently regular and can be directly computed from the coefficients defining the Hermite interpolant. In our tests we consider $\dt$ small enough such that this elastic regularization is warranted. Extensions for larger time step predictions incorporating numerical methods for the flows of diffeomorphisms \cite{ashburner2011diffeomorphic, miller2006geodesic, beg2005computing} and a more exhaustive numerical study for other energy functionals has been deferred to future investigation. Our numerical experiments presented here can be seen as proof of concept for the use of the lifting operator approximation and the functional discretization \eqref{diffeomorphic_functional_discretization}. \par

We implement the lifting operator approximation \eqref{neural_lifting_operator} using the FNO architecture with the open-source Neural Operator package\footnote{\href{https://github.com/neuraloperator/neuraloperator}{https://github.com/neuraloperator/neuraloperator}}. The formulation however enables the application of other existing approximations \cite{lu2021learning, bhattacharya2021model, raonic2023convolutional, tripura2022wavelet, liu2024neural} with suitable modification. As comparison, we consider the straightforward application of the FNO architecture without the lifting operator construction, referred to as the FNO-fields solution and our implemented approach as the FNO-maps solution. The FNO-fields solution is trained using only the discretized similarity term in the cost function. Each architecture takes input data with five stacked input fields, uses a Fourier wavenumber truncation at $k = 32$ in each coordinate direction with four Fourier layers, and $128$ hidden channels. The FNO-maps architecture uses the same parameters with eight output channels for the components defining the interpolant \eqref{diff_interpolant}. We have included an open-source repository\footnote{Code made available upon acceptance: \href{https://github.com/seth-tl/Diff_flow_learning}{https://github.com/seth-tl/Diff\_flow\_learning}} containing all the code used to implement these experiments. The optimization is performed using the built in PyTorch Adam optimizer \cite{da2014method} with a learning rate of $10^{-3}$ and without significant hyper-parameter tuning. The data was generated using the CMM for advection on the torus \cite{mercier2020characteristic} and for Euler's equations \cite{yin2021characteristic} and the number of sample points and resolution used in training is specified for each problem.\par 

We perform the time stepping with the scheme \eqref{composition_time_stepping} using the lifting operator approximation \eqref{neural_lifting_operator} which requires an interpolation of the initial field $u_0$ and we consider a bilinear interpolant. The composition required to perform the time stepping is evaluated using the Hermite interpolant representation of the maps. We found that the associated $\mathcal{O}(N_t \cdot |D_h|)$ memory footprint in storing the history of the maps was negligible for the number of time steps considered. 

The error is assessed with the mean squared error against both in and out of distribution samples. We additionally measure the conservation error
\begin{equation}
\label{conservation_integral}
\text{conservation error }(\tilde{a}(t)) =  \int_{\mathbb{T}^2} (a_0\circ \tilde{\varphi}_{[t,0]} \cdot \text{det}(D\tilde{\varphi}_{[t,0]}) - a_0) dx \,,
\end{equation}
where $\tilde{\varphi}_{[t,0]}$ is the learned solution map.  Due to the functional nature of the evolution, the method is conservative in a continuous sense and the integral \eqref{conservation_integral} must be evaluated via quadrature. We consider a quadrature scheme based on Simpson's rule on a fine grid however the errors cannot be absolved entirely of this quadrature approximation in the presence of discontinuities or for the highly multi-scale fields output by the method. 

\subsection{Transport of a discontinuous field}

\begin{figure}[h!]
\centering
\includegraphics[width = 8cm, height = 3cm]{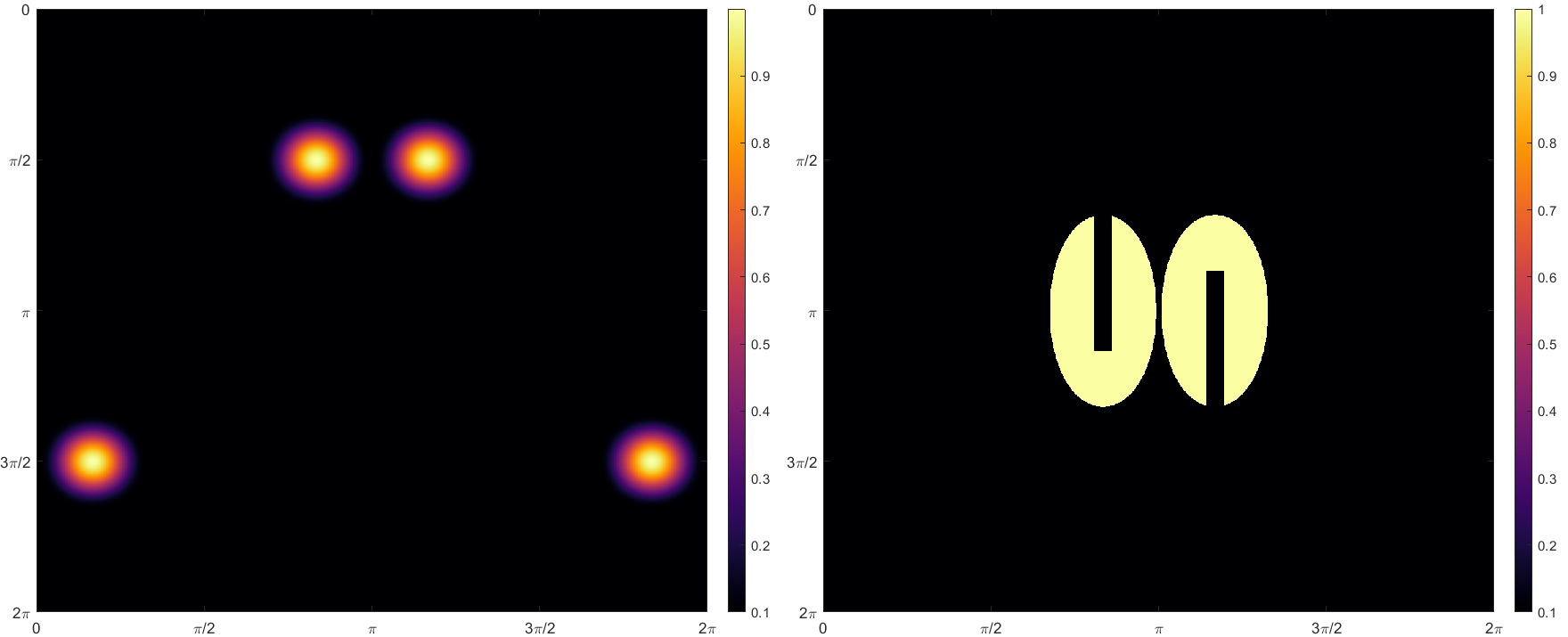}
\caption{Training field (left) and testing field (right) for the linear advection test case.}
\label{linear_advection_initial_conditions}
\end{figure}

\begin{figure}[h!]
\centering
\includegraphics[width = 10cm, height = 4.5cm]{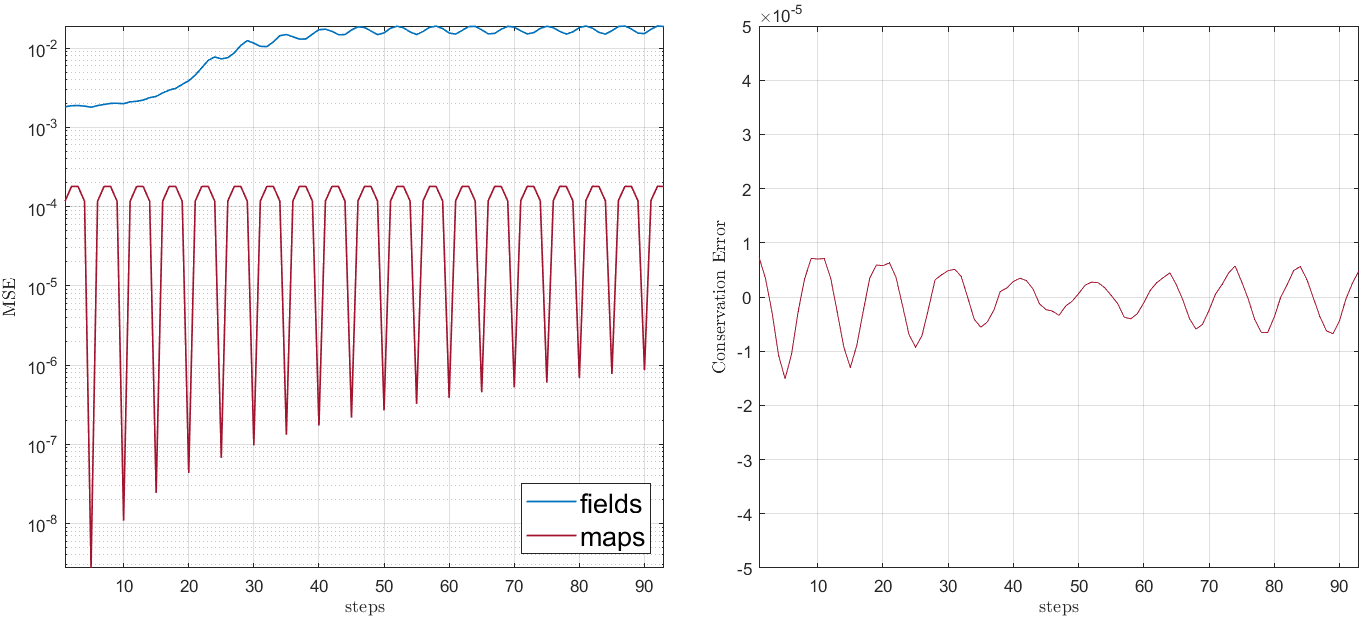}
\caption{Left: Error over time for the linear advection out of distribution test.  Right: Conservation of mass error over time}
\label{fig:errors_over_time_advection}
\end{figure}

\begin{figure}[h!]
\centering
\includegraphics[width = \linewidth, height = 8cm]{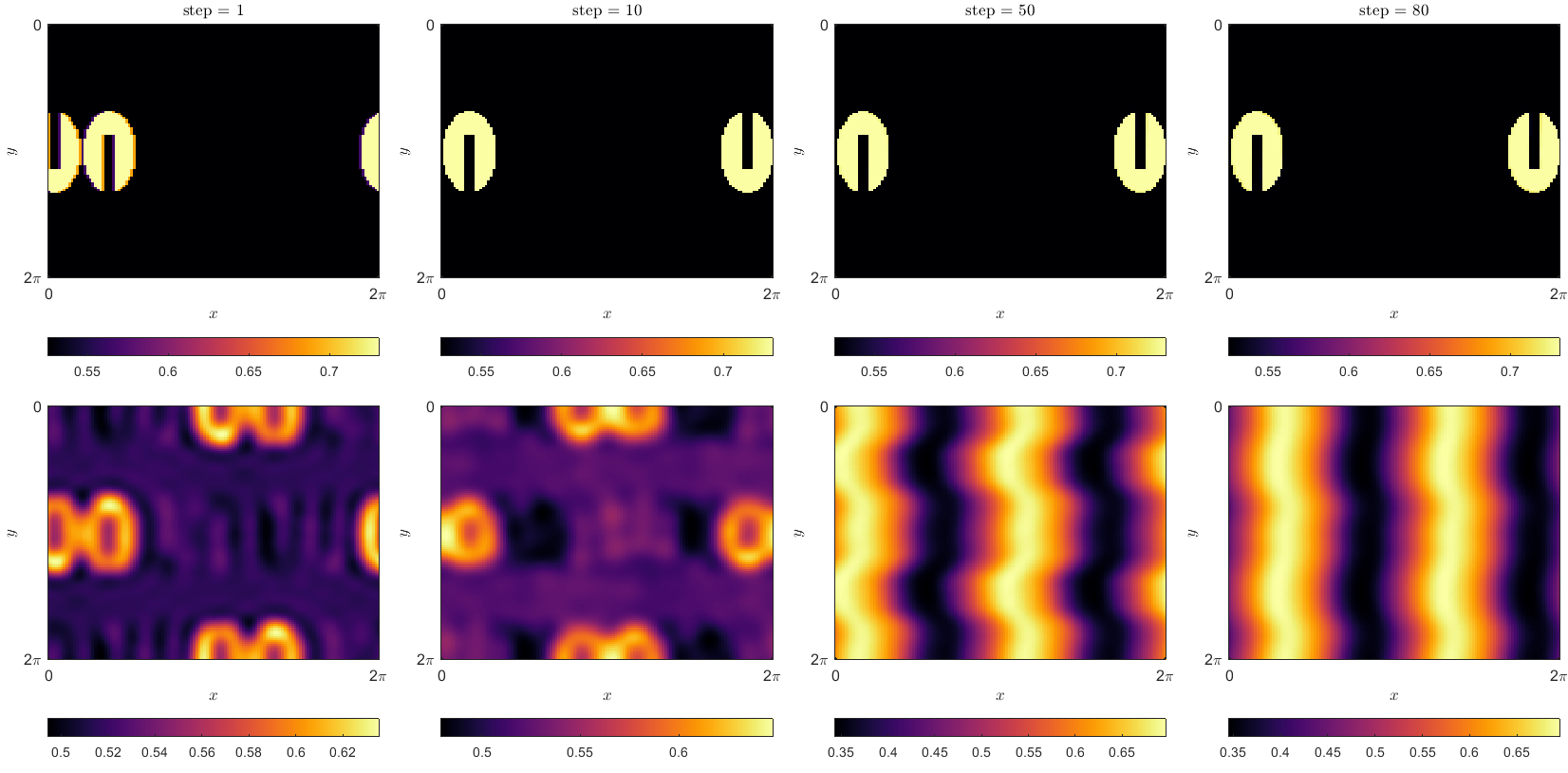}
\caption{Linear advection test case, evolution of the solution. Top row: Solution from learned lifting operator learning. Bottom row: Solution learning the field evolution directly.}
\label{fig:time_evolution_errors_advection}
\end{figure}

As a first test we consider learning the solution operator to the linear advection equation
\begin{equation}
\partial_t \phi + \bsym{u} \cdot \nabla \phi = 0 \,, \quad \phi(0) = \phi_0 \,,
\end{equation}
in one of the simplest possible scenarios with the velocity field defining a uniform translation $\bsym{u} = [2\pi,0]^T$. We train on a single initial condition and assess the ability of the network to generalize to out-of-distribution fields by transporting a discontinuous slotted cylinder (see Figure ~\ref{linear_advection_initial_conditions}). The ground truth data is obtained using a single simulation sampled at time steps $\Delta t = 0.01$ over a time interval $[0,10]$ on a grid of size $128 \times 128$. The test is meant to demonstrate how the method rapidly generalizes to out-of-distribution transported quantities by learning the transformation of the transported field instead of predicting the changing intensity values. In this scenario, the lifting operator is simply the constant $u \mapsto [ 2\pi \dt,0]^T$. The FNO-fields and FNO-maps solutions were able to achieve machine-precision in the loss functional at training time and in Figure~\ref{fig:time_evolution_errors_advection} we illustrate the evolution of the error in the discontinuous initial condition applied to learning the lifting operator and learning the evolution operator directly. Using our proposed approach we see that the evolution is non-diffusive and retains the discontinuity in the initial condition. The right panel of Figure~\ref{fig:errors_over_time_advection} shows the conservation of mass errors, observing an approximate conservation of mass for the discontinuous initial condition. 

\subsection{Turbulent energy cascades}

\begin{figure}[h!]
\centering
\includegraphics[width = \linewidth, height = 6cm]{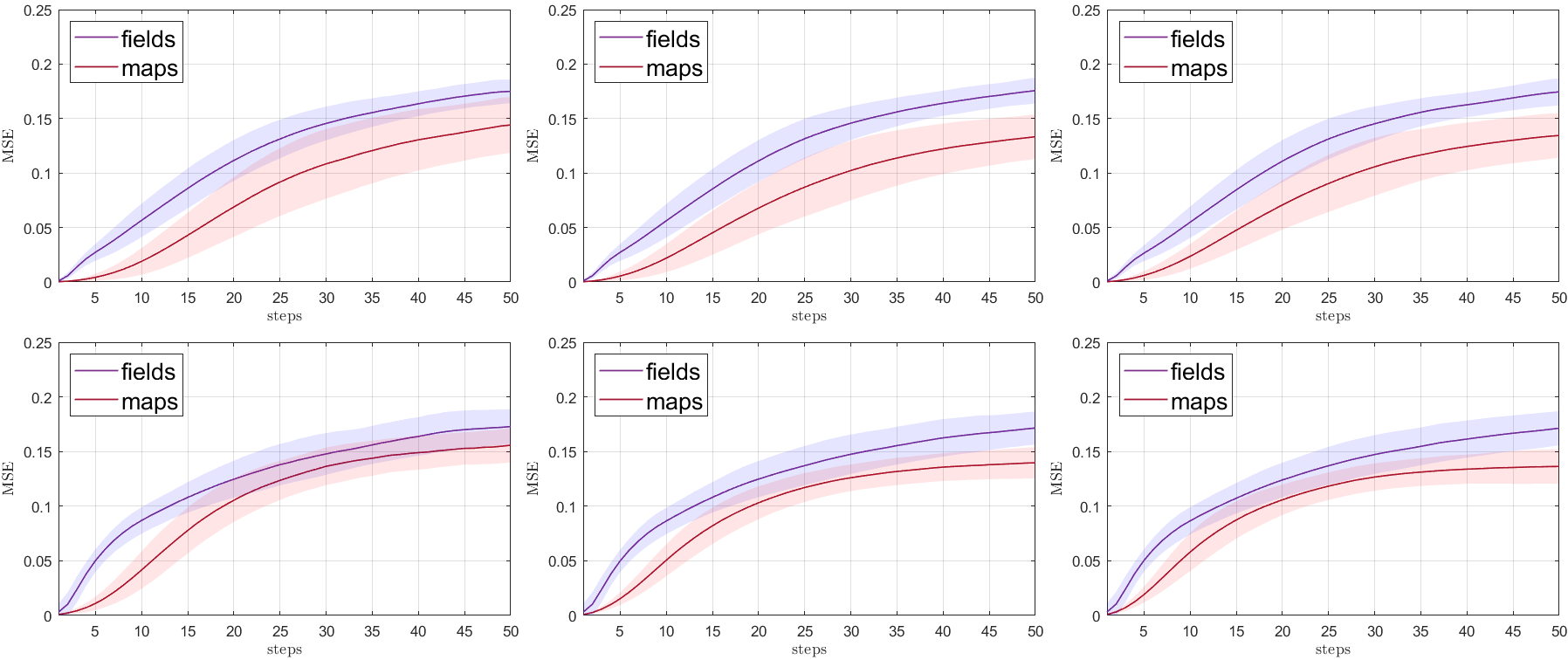}
\caption{Mean and standard deviation of the error over time for in-distribution (top row) and out-of-distribution samples (bottom row). Time-stepping applied auto-regressively for the trained operator networks with inference performed on fields with spatial resolution $\Delta x = 2\pi / 128, 2\pi / 256, 2\pi / 512$ from left to right.}
\label{fig:errors_euler}
\end{figure}

\begin{figure}[h!]
\centering
\includegraphics[width = 8cm, height = 5cm]{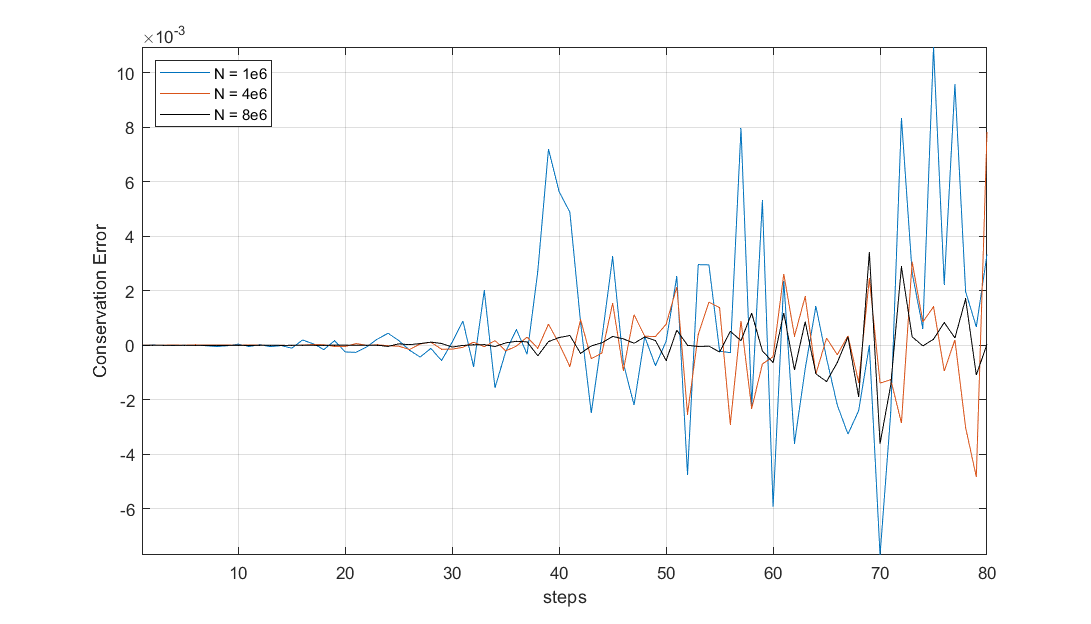}
\caption{Enstrophy conservation error over time evaluated a different quadrature grids.}
\label{fig:conservation_errors_Euler}
\end{figure}

In this test case we consider learning the solution operator to the two-dimensional incompressible Euler equations 
\begin{equation}
\label{euler_eqs}
\partial_t \bsym{u} + \bsym{u}\cdot \nabla \bsym{u} = -\nabla p \,,\quad \nabla \cdot \bsym{u} = 0 \,, \quad u(0) = \bsym{u}_0 \,.
\end{equation}
These equations are a canonical example of turbulent fluid dynamics exhibiting a forward and inverse cascade of energy and are notoriously difficult to accurately solve numerically despite being the subject of decades of study. Many of the computational difficulties associated with the equations come from the rapid generation of vorticity gradient which create fine scale features which can exceed the resolution of the computational grid. The vorticity field $\omega(t) = \nabla^{\perp}\bsym{u}$ is a transported quantity, evolving as $\omega(t) = \omega_0 \circ \varphi_{[t,0]}$ and we consider learning the solution operator based on sequence of vorticity fields. The training data was generated using initial vorticity distributions of the form
\begin{equation}
\label{initial_vorticity}
\omega_0(x) = \sum_{|k| \leq K} a_k\cos(k\cdot x) + b_k \sin(k\cdot x) \,,
\end{equation}
with the constants $a_k$ and  $b_k$ drawn from a uniform distribution over $[-1,1]$. The simulations were performed up to a final integration time of $T = 10$ with $\Delta t = 0.001$, a spatial resolution $\Delta x = 2\pi/128$ and remapping every 10 time steps with $\Delta \tau = 10 \dt$ which defines the time discretization in the dataset \eqref{dataset} samples. We consider an 80\% / 20\% train / test data split with initial conditions defined such that $K \leq 10$, over a total of $100$ simulations and trained each network over $50$ epochs. We additionally tested the networks against out-of-distribution for another $20$ simulations defined by initial conditions with $K \leq 20$. The errors of the method applied to both test data sets are shown in Figure~\eqref{fig:errors_euler}.  The error saturates for both in-distribution and out-of-distribution tests after approximately 25 steps and we observe that the method outperforms the solution trained on the field data alone in both cases. We tested the super-resolution properties of each architecture, evaluating the errors over time at inference on vorticity samples at two and four times finer resolution. The error is observed to behave independent of the resolution, demonstrating super-resolution properties of the neural operators \cite{li2020fourier}. In figure \ref{fig:conservation_errors_Euler} we demonstrate the conservation error as evaluated by \eqref{conservation_integral} for the square vorticity which demonstrates an approximate analytic conservation over time. \par

\begin{figure}[h!]
\centering
\includegraphics[width = 14cm, height = 10cm]{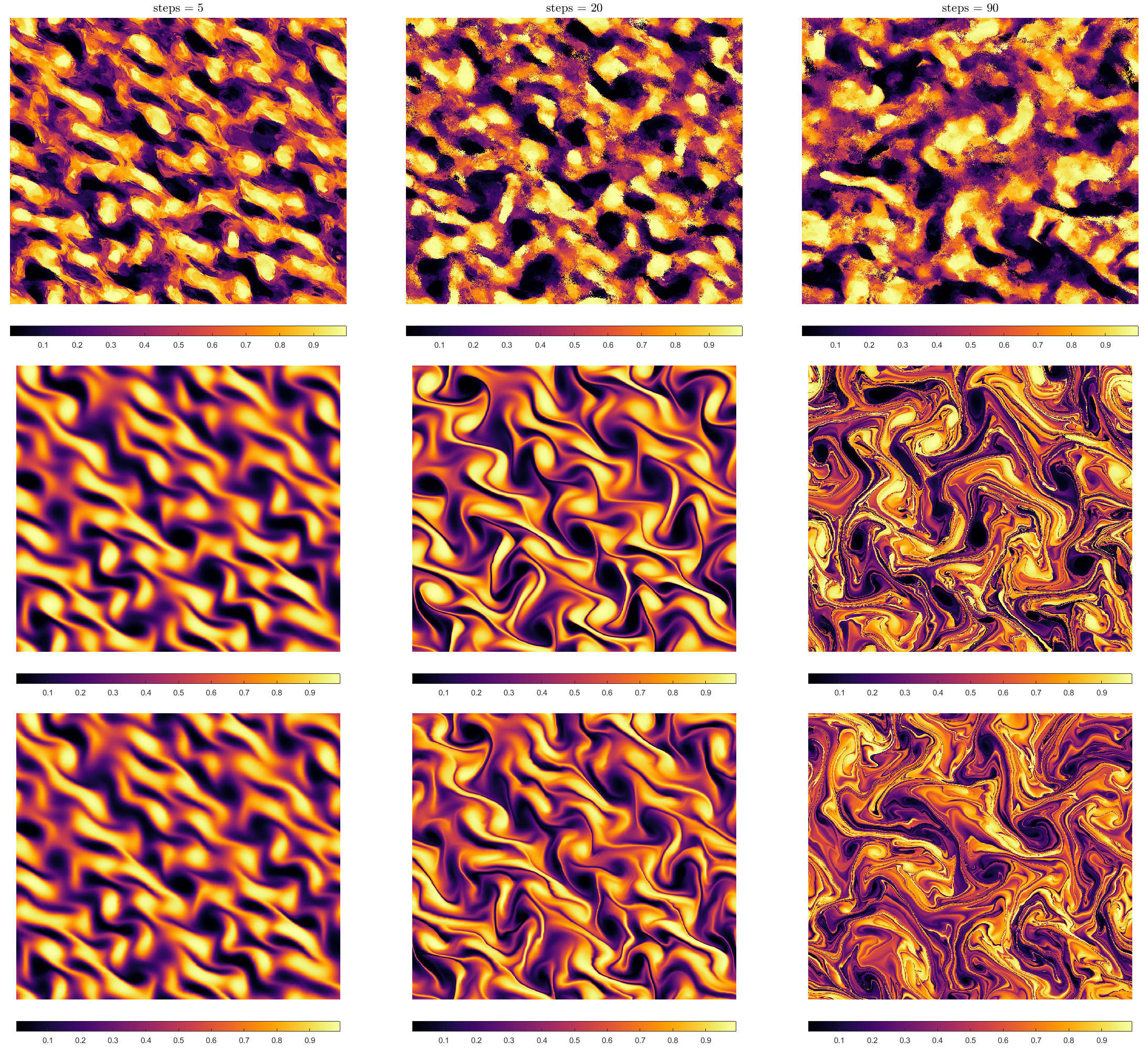}
\caption{Evolution of test solution at three snapshots in time. Top row: FNO-fields solution. Middle: FNO-maps solution. Bottom row: Ground truth. The statistical properties of the solution are seen to be retained with our learned solution along with fine scale vortex structures over many time steps. }
\label{fig:turbulent_evolution}
\end{figure}

\begin{figure}
\centering
\includegraphics[width = \linewidth, height = 6cm]{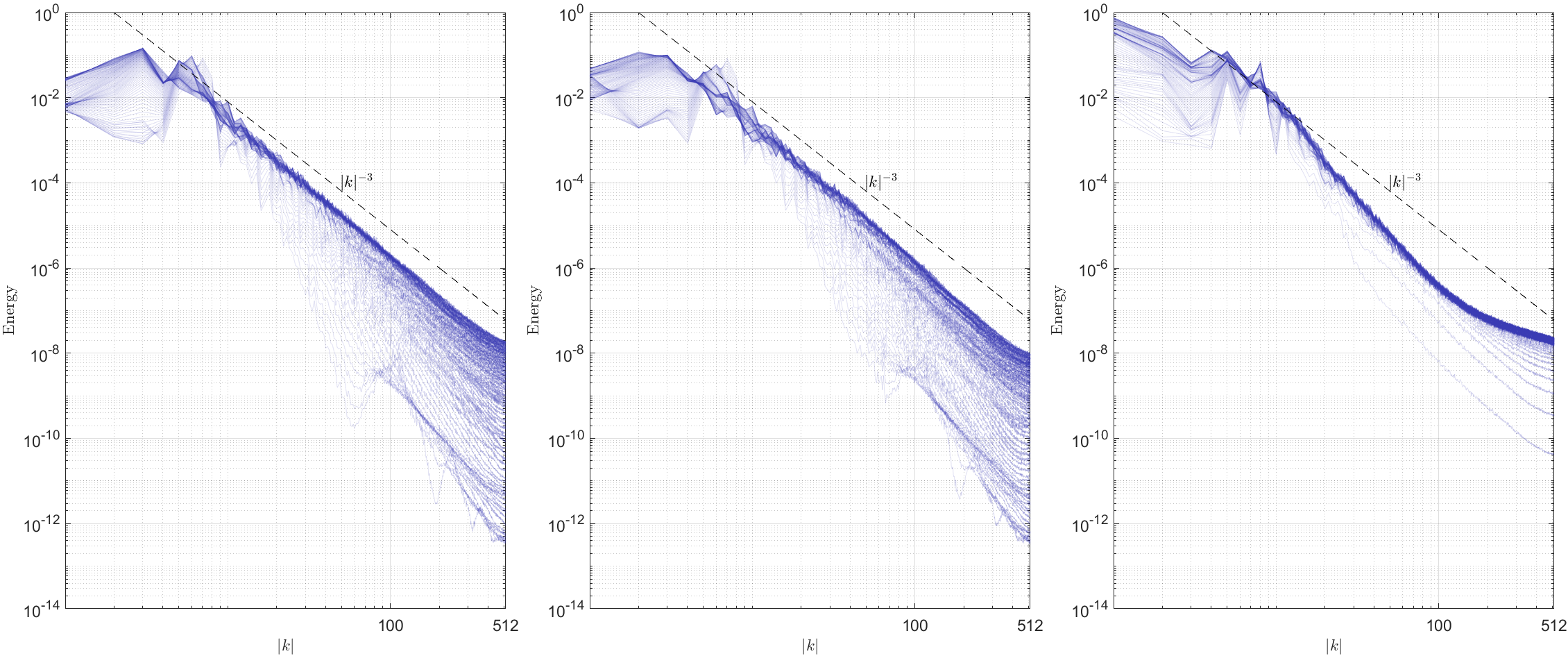}
\caption{Evolution of the energy spectrum for the decaying turbulence experiment associated with the flow \ref{fig:turbulent_evolution}. Left: Ground Truth. Middle: FNO-maps solution. Right: FNO-fields solution.}
\label{fig:energy_cascades}
\end{figure}

A quantitative evaluation of the super-resolution properties of each neural operator architecture can be given by measuring the energy cascades at higher-resolution. Turbulent fluid dynamics possess a self-similarity across spatial scales \cite{kraichnan1980two} resulting in a well-defined power law in the energy spectrum. In the case of two-dimensional incompressible fluid dynamics, a $|k|^{-3}$ scaling is expected to form in the energy spectrum as the fluid transitions into a turbulent regime \cite{boffetta2012two}. We demonstrate the ability of our proposed strategy to reproduce the expected cascade scaling at subgrid scales by performing the inference on a resolution of $1024\times 1024$ for the fields solution and by upsampling the output maps for an inference on the same training resolution. The evolution is visualized in Figure~\ref{fig:turbulent_evolution} where we observe stable energy cascade behaviour in our solution for our proposed approach. We note that the presented dataset possesses fields with a large range of scales and the FNO-fields network is observed to have learned both noisy features along with the large scale dynamics. In Figure~\ref{fig:energy_cascades} we plot the energy spectrum of the solution over time for a single sample. The upsampled numerical solution obtained via our proposed learning strategy is observed to capture the expected cascade scaling at subgrid scales. To the best of our knowledge, no other operator learning method possesses these resolution properties.  

\section{Conclusion}
\label{sec:conclusions}

In this work we have formulated an operator learning strategy for evolution equations based on a lifting into the space of diffeomorphisms.  The operator learning problem \eqref{problem1} is transformed into an auxiliary problem by approximating the evolution operator of a dynamical system in a factorized form with a learned lifting operator and the action of the diffeomorphism group on the field space. This technique permits the direct application of existing network architectures while enforcing a number of desirable properties into the learning problem as hard constraints. In particular, we have shown that the approximation enforces provable conservation principles and allows for the approximations of turbulent fluid flows with the correct statistical scaling laws in the energy spectrum. We have validated the methodology on some canonical test cases, observing state of the art performance for evolution operator learning of turbulent fluid flows.

\textbf{Extensions.} This work has focused on learning evolution operators which contain (nonlinear) transported quantities without the presence of diffusion or forcing. The methodology however admits some immediate extensions. A simple residual based technique can be used to model the dynamics which are not transported. In particular, the evolution operator can be decomposed as 
\begin{equation}
\Phi_t  = \Psi \circ (\mathcal{K}_t, I) + \mathcal{A}_t
\end{equation}
where $I: U \to U$ is the identity and $\mathcal{A}_t : U \to U$ is a residual operator which is also learned. Another decomposition could take the form of a splitting-based approach to model the effects of diffusivity
\begin{equation}
\Phi_t = \mathcal{D}_{t} \circ \Psi \circ (\mathcal{K}_t, \mathcal{D}_{t}) \,,
\end{equation}
where $\mathcal{D}_{t}:U \to U$ capture the dynamics in the Lagrangian frame of reference. The inclusion of source terms can be made explicitly through Duhamel's principle. In particular, consider a generic advection equation of the form 
\begin{equation}
\partial_t \phi + \bsym{u} \cdot \nabla \phi = f(t)
\end{equation}
where the source term may have nonlinear dependence on the transported quantity. It can be shown that the evolution of $\phi$ can be decomposed as 
\begin{equation}
\phi(t) = \phi_0 \circ \varphi_{[t,0]} + \varphi_{[t,0]}^* \int_0^{t} f(s) \circ \varphi_{[0,s]} ds = \left (\phi_0 + \int_0^t f(s) \circ \varphi_{[0,s]} ds\right) \circ \varphi_{[t,0]}
\end{equation}
The accumulated source term operator
\begin{equation}
F_{[s,t]}: U \to U \,, \quad \phi \mapsto \int_s^t f(s, \phi(s)) \circ \varphi_{[0,s]} ds \,,
\end{equation}
could be incorporated in the learning problem such that the evolution operator is decomposed as 
\begin{equation}
\Phi_t = \Psi \circ (\mathcal{K}_t, I + F_{[0,t]}) \,.
\end{equation}
Since the source term operator can be decomposed additively, the resulting time-stepping scheme becomes
\begin{equation}
u_{n+1} = \left(u_0 \circ \mathcal{K}_{1}(u_0) \circ \dots \circ \mathcal{K}_{{n-1}}(u_{n-1}) + F_{[0,n-1]}(u_{n-1}) + F_{[n-1,n]}(u_n) \right) \circ \mathcal{K}_{n}(u_n) \,.
\end{equation}
As indicated in section \ref{sec:LDM_connection}, the operator learning strategy we proposed directly connects to the ideas of shape analysis \cite{younes2010shapes} and diffeomorphic registration. In a similar spirit to the work of Benn and Marsland \cite{benn2022measurement}, we hope that further insights can be leveraged from shape analysis for the design of operator learning strategies. Multi-scale techniques \cite{bruveris2012mixture} which consider fine and coarse scale diffeomorphisms in the context of downscaling along with the techniques of metamorphosis \cite{holm2009euler} to incorporate changing topologies could be fruitful lines of investigation.

\textbf{Outlook.} Our methodology demonstrates an example of how purely data-driven approaches to evolution operator learning problem can be transformed to enforce known physical laws directly. This technique can be view as a form of \emph{geometric operator learning}, building off the approaches of geometric deep learning \cite{bronstein2017geometric} to an infinite-dimensional setting. We have considered a particular geometric operator learning strategy to enforce a relabelling symmetry into the problem and it is our hope that similar approaches can be developed to embed other invariants, such as time-translation for the conservation of energy, directly into the functional form of the operator. Given the burgeoning applications of deep learning in operational weather forecasting, advancing the design of operator learning techniques which can directly incorporate known physical laws as hard constraints should enable new hybrid modelling techniques. Assessing the benefits of this approach through its extension for nonlinear dynamic forecasting by training on ERA5 re-analysis data \cite{hersbach2020era5} would bring further insights into its limitations and capabilities for numerical weather prediction and climate modelling. 

\section*{Acknowledgments}
This research was undertaken, in part, thanks to funding from the Canada Research Chairs program, the NSERC Discovery Grant program, and the NSERC CGS-D program. ST would like to thank Jonathan Chalaturnyk for the helpful suggestions and insights on the implementation.

\FloatBarrier

\small
\bibliographystyle{unsrt}  
\bibliography{Bibliography}

\appendix

\section{Proofs for Section \eqref{sec:method}}
\label{app:proofs}

\subsection{Section \ref{sec:resolution_properties}}
\begin{proof}(of Lemma \ref{lem:composite_bandwidth}). 
Consider the composition of two elements $\varphi_1, \varphi_2 \in D_h$ and write $\varphi_1 = \id + u$ and $\varphi_2 = x + v$ such that $\varphi_1 \circ \varphi_2 = \id + v + u\circ(\id + v)$. Since $u \in \mathcal{B}_L^d$ we can determine the $\epsilon$-effective bandwidth of the last term by writing 
\begin{equation}
\label{expanded_composition}
u \circ (x + v(x))  = \sum_{|k| \leq L} \hat{u}_k e^{ik \cdot x}\sum_{n = 0}^{R-1} \frac{(ik \cdot v(x))^n}{n!} + T_R(x) \,,
\end{equation}
and computing $R$ such that
\begin{equation}
\label{tail_bound}
\norm{T_{R}}_2  = \bigg\|  \sum_{|k| \leq L} \hat{u}_k e^{ik \cdot x} \sum_{n = R}^{\infty} \frac{(ik \cdot v(x))^n}{n!} \bigg\|_2  \leq \norm{u}_{L^2} \frac{(L \norm{v}_{\infty})^{R}}{R!} \frac{1}{1- \frac{L\norm{v}_{\infty}}{R + 1}} \leq \epsilon \,,
\end{equation}
with the condition that $R + 1 > L \norm{v}_{\infty}$. Using the upper bound $ n! < \sqrt{2\pi n}n^n e^{-n}e^{1/12n } $ valid for all integers \cite{robbins1955remark}, we get that $R \geq L \norm{v}_{\infty} + \log(\norm{u}_2/\epsilon)$ is sufficient to satisfy the inequality (neglecting the $1/12n$ correction). Then since $v \in \mathcal{B}^d_L$ it follows that the first term in \eqref{expanded_composition} is in $\mathcal{B}^d_{RL + L}$ which gives 
\begin{equation}
\text{bw}_{\epsilon}(\varphi_1 \circ \varphi_2 - \id) \leq  L + L\left(L \norm{v} _{\infty} + \log( \norm{u}_2/\epsilon)\right)
\end{equation}
Using this same reasoning, the bound \eqref{composite_bandwidth} for the $k$-times composition satisfies 
\begin{equation}
\text{bw}_{\epsilon}(\varphi_1 \circ \cdots\circ\varphi_k - \id) \leq \text{bw}_{\epsilon}(\varphi_2 \circ \cdots\circ\varphi_k - \id)(1 + L\norm{\varphi_2 \circ \cdots \circ \varphi_k - \id }_{\infty})  + L \log(\norm{\varphi_1 - \id}_2/\epsilon) + L  \,.
\end{equation}
Let $C_j = \norm{\varphi_j - \id}_2$ and $B_j = \norm{\varphi_{j-1} \circ \cdots \circ \varphi_k - \id}_{\infty}$, then we can iterate the recurrence relation to give 

\begin{equation}
\text{bw}_{\epsilon}(\varphi_1 \circ \cdots \circ \varphi_k - \id) \leq L \sum_{i = 1}^k \left(1 + \log(C_i/\epsilon)\right) \prod_{j = 2}^i(1 + LB_j)
\end{equation}
with the $i =1$ product is set to one. Uniformly bounding the constants $C_j$ and $B_j$ yields the bound \eqref{composite_bandwidth}. The bound \eqref{orbit_space_estimate} follows as corollary using a similar computation.
\end{proof}

\begin{proof}(of Lemma \ref{lem:resolution_consistency_Diff}.) Let $\{M_L\}_{L = 1}^{\infty}$ be a discrete refinement and define a reconstruction operator $\mathscr{R}_L: \mathbb{R}^{L \cdot m} \to U$ which satisfies the properties: (P1) $\mathscr{S}_L \circ \mathscr{R}_L \circ \mathscr{S}_L = \mathscr{S}_L$ and the consistency condition (P2) $\lim_{L \to \infty}\norm{\mathscr{R}_L \circ \mathscr{S}_L(u) - u}_U = 0$ for all $u \in K \subset U$ compact. Define the discrete action $\tilde{\Phi}_L(M_L, \mathscr{S}_L(u), \theta) = \Psi(F(\theta), \mathscr{R}_L \circ \mathscr{S}_L(u))$ and note that the discretization invariance \eqref{discretization_invariance} with the continuum operator $\Psi(F(\theta), \cdot)$ holds immediately due to continuity of the group action and (P2). The resolution consistency then holds since
\begin{equation}
\begin{aligned}
\mathscr{S}_{L'}\left(\tilde{\Phi}_L(M_L,\mathscr{S}_L(u), \theta)\right) &= \mathscr{S}_L'\left(\Psi(F(\theta), \mathscr{R}_L \circ \mathscr{S}_L(u)) \right) =  \mathscr{S}_L'\left(\Psi(F(\theta), \mathscr{R}_{L'} \circ \mathscr{S}_{L'}(u)) \right) 
\\
&= \mathscr{S}_{L'}\left( \tilde{\Phi}_{L'}(M_{L'}, \mathscr{S}_{L'}(u), \theta) \right)
\end{aligned}
\end{equation}
where the second equality uses the fact that the group action is pointwise defined and the property (P1).  
\end{proof}

\subsection{Section \ref{sec:LDM_connection}}

\begin{proof}(of Lemma \ref{lem:measurable_selector}.)
The claim follows by application of the Kuratowski–Ryll–Nardzewski measurable selection theorem (see Theorem 18.13 \cite{aliprantis2006infinite}). In the particular form needed here, it states that: for $X$ measurable and $Y$ a complete and separable metric space, then a correspondence $F: X\rightrightarrows Y$ with non-empty, closed values and a measurable graph admits a measurable selector $f: X \to Y$ such that $f(x) \in F(x)$ for all $x \in X$. In our case, $X = U$ is a Hilbert space, and $Y = \Diffs$ is a complete and separable metric space \cite{bruveris2017completeness}, so we must demonstrate closedness of the values of the lifting operator \eqref{lifting_operator_definition}. Let $(u_n) \subset K \subset H^{s+1}$ be a bounded sequence converging to $u$ and consider some $\varphi_n \in \mathcal{K}_{\dt}(u_n)$. Due to the definition \eqref{lifting_operator_definition} it follows that 
\begin{equation*}
\frac{1}{2}\dist^2_{s}(\id, \varphi_n) \leq  \frac{1}{2}\dist^2_{s}(\id, \varphi_n) + \norm{\Psi(\varphi_n,u_n) - \Phi_{\dt}(u_n)}^2_{s} \leq   \norm{u_n - \Phi_{\dt}(u_n)}^2_{s} \leq C\,,
\end{equation*}
for some uniform constant depending on the time step size. Using Theorem 7.22 of \cite{younes2010shapes}, there exists a $\bv_n \in L^2([0,T],H^s) = X^{s,2}([0,T])$ with $\varphi_n$ being its time-$T$ flow such that 
\begin{equation}
\label{v_hilbert_bound}
\dist^2_s(\id, \varphi_n) = \int_0^T\norm{\bv_n(t)}_{s}^2 dt. 
\end{equation}
Since $X^{s,2}([0,T])$ is a Hilbert space and the sequence $(\bv_n) \subset X^{s,2}([0,T])$ is uniformly bounded, there exists a weakly convergent subsequence $\bv_{n_k} \rightharpoonup \bv \in X^{s,2}([0,T])$ which also weakly converges in the larger space $X^{s,1}([0,T])$. This implies strong convergence in $X^{s-1,1}([0,T])$ by compactness of the embedding $H^{s} \hookrightarrow H^{s-1}$. This implies strong convergence of the flows $\varphi_n \to \varphi$ in $\text{Diff}^{s-1}(M)$ and $\varphi \in \Diffs$ by Lemma 4.2 \cite{bruveris2017completeness}. Using lower semi-continuity of the distance and continuity of the evolution operator and group action, along with the fact that each $\varphi_n$ is a minimizer gives us that $\varphi$ satisfies
\begin{equation*}
\begin{aligned}
\frac{1}{2}\dist^2_{s}(\id, \varphi) + \norm{\Psi(\varphi,u) -\Phi_{\dt}(u)}_{s}^2 &\leq \liminf_{n \to \infty} \frac{1}{2}\dist^2_{s}(\id, \varphi_n) + \norm{\Psi(\varphi_n,u_n) - \Phi_{\dt}(u_n)}_s^2
\\
&\leq \frac{1}{2}\dist^2_{s}(\id, \eta) + \lim_{n \to \infty} \norm{\Psi(\eta,u_n) - \Phi_{\dt}(u_n)}_s^2 
\\
&= \frac{1}{2}\dist^2_{s}(\id, \eta) + \norm{\Psi(\eta,u) - \Phi_{\dt}(u)}_s^2 \,, \quad \text{ for all  } \eta \in \Diffs \,,
\end{aligned}
\end{equation*}
and therefore $\varphi \in \mathcal{K}_{\dt}(v)$ which establishes the claim. 
\end{proof}

\begin{proof}(of Theorem \ref{thm:approx_thm1}.) Note that the bound \eqref{set_control_Diff} gives us that the measurable selector is also square-integrable. Define $\mathcal{V}: U \to \mathcal{U} \subset \mathfrak{X}^s(M)$ such that $\mathcal{V}(u) = J^{-1} \circ \tilde{\mathcal{K}}_{\dt}(u)$. Since $J^{-1}$ is Lipschitz, we have that $\mathcal{V}$ is $\mu$-square integrable and using the existence Theorem 18 of Kovachki et al. \cite{kovachki2021universal} (see also Lanthaler et al. \cite{lanthaler2304nonlocality}) we get that there exists a neural operator $\mathcal{N}_{\theta}: U \to \mathcal{U} \subset \mathfrak{X}^s(M)$ such that 

\begin{equation}
\label{error_L2_NO}
\int_U \norm{\mathcal{V}(u) - \mathcal{N}_{\theta}(u)}_{s}^2 d\mu(u) < \varepsilon'
\end{equation}
for any $\varepsilon' > 0$. Then note that the first term in \eqref{operator_approx_error} is introduced by splitting the error with the selector $\tilde{\mathcal{K}}_{\dt}: U \to \Diffs$ and we can use the differentiability of the composition operator \cite{inci2013regularity} and Lipschitz continuity of the chart function to give 
\begin{equation*}
\begin{aligned}
\norm{ u \circ \tilde{\mathcal{K}}_{\dt}(u) -  u \circ  J \circ \mathcal{N}_{\theta}(u)}_{s} 
&\leq C(\norm{u}_{s+1})\dist_{s}(\tilde{\mathcal{K}}_{\dt}(u), J \circ \mathcal{N}_{\theta}(u))
\\
& \leq  C(\norm{u}_{s+1}, L_{J})\norm{\mathcal{V}(u) - \mathcal{N}_{\theta}(u)}_s \,.
\end{aligned}
\end{equation*}
Given $\varepsilon> 0$, using \eqref{error_L2_NO}, and restricting to $K \subseteq B_R(0) \subset U$, we can take $\varepsilon' < \varepsilon/C(R, L_J)$ to establish the claim.
\end{proof}

\begin{proof}(of Corollary \ref{orbit_time_stepping_corollary}.) Using the conditions imposed on the evolution operator, we can can construct a continuous representative \eqref{continuous_rep_lift}, which satisfies the factorization \eqref{evolution_operator_decomposition} exactly by considering the constrained optimization problem 
\begin{equation}
\mathcal{K}_{\dt}(u) = \arg\min_{\varphi \in \Diffs} \frac{1}{2}\dist^2_s(\id ,\varphi) \,, \quad \text{s.t.} \quad u \circ \varphi = \Phi_{\dt}(u) \,,
\end{equation}
and an associated measurable selector $\tilde{\mathcal{K}}_{\dt}(u) \in \mathcal{K}_{\dt}(u)$ such that $\Psi(\tilde{\mathcal{K}}_{\dt}(u),u) = \Phi_{\dt}(u)$. The proof proceeds similarly to the proof of theorem \ref{thm:approx_thm1}. Let $u_{n+1} = u_{n} \circ \tilde{\mathcal{K}}_{\dt}(u_{n})$ and $\tilde{u}_{n+1} = u_{n} \circ J \circ \mathcal{N}_{\theta}(\tilde{u}_{n})$ both initialized as $\tilde{u}_0 = u_0 = u$. Using the Sobolev embedding the $H^s$ norm is sufficient to control the $L^{\infty}$ norm with $s > d/2$. We can write the pointwise bound for the error $u \in \bar{K}$ as

\begin{equation*}
\begin{aligned}
\norm{\Phi_{\dt}^{\circ,n+1}(u) - \Psi \circ J \circ \mathcal{N}_{\theta}^{\circ,n+1}(u)}_{\infty}
& = \norm{u_{n} \circ \tilde{\mathcal{K}}_{\dt}(u_{n}) - \tilde{u}_{n} \circ J \circ \mathcal{N}_{\theta}(\tilde{u}_{n})}_{\infty}
\\
&\leq C\norm{u_n}_{s} \norm{\tilde{\mathcal{K}}_{\dt}(u_n) - J \circ \mathcal{N}_{\theta}(\tilde{u}_n)}_{\infty} + \norm{u_n - \tilde{u}_n}_{\infty}
\end{aligned}
\end{equation*}
We can then split the first term as 
\begin{equation*}
\begin{aligned}
\norm{\tilde{\mathcal{K}}_{\dt}(u_n) - J \circ \mathcal{N}_{\theta}(\tilde{u}_n)}_{\infty} &\leq L_J(\norm{\mathcal{V}(u_n) - J \circ \mathcal{N}_{\theta}(u_n)}_{\infty} + \norm{\mathcal{N}_{\theta}(\tilde{u}_k) - \mathcal{N}_{\theta}(u_k)}_{\infty})
\\
&\leq  L_J(\norm{\mathcal{V}(u_n) - J \circ \mathcal{N}_{\theta}(u_n)}_{\infty} + L_{\mathcal{N}}\norm{\tilde{u}_k - u_k}_{\infty})
\end{aligned}
\end{equation*}
where $L_{\mathcal{N}}(\theta)$ is the Lipschitz constant of the neural operator over $\bar{K}$. Iterating with the previous bound gives 
\begin{equation}
\norm{u_{n+1} - \tilde{u}_{n+1}}_{\infty} \leq C(R)\sum_{j = 0}^n (1 + L_JL_{\mathcal{N}})^{n-j}\norm{\mathcal{V}(u_k) - \mathcal{N}_{\theta}(u_k)}_{\infty} \,.
\end{equation}
Then using the bound on $\bar{K}$ we can invoking \eqref{error_L2_NO} again with $\varepsilon = \varepsilon'$ and use linearity of the integration to establish the claim.  
\end{proof}

\end{document}